\documentclass[a4paper,11pt]{article}
\usepackage{mathrsfs}
\usepackage{latexsym}
\usepackage{amsmath,amssymb}
\usepackage[pdftex]{hyperref}
\usepackage{amsthm}
\usepackage{amsfonts,cite}
\usepackage[usenames]{color}
\usepackage{amssymb}
\usepackage{graphicx}
\usepackage{amsmath}
\usepackage{amsfonts}
\usepackage{amsthm}
\usepackage{mathrsfs}
\usepackage{dsfont}
\usepackage{indentfirst}

\newtheoremstyle{mythm}{1.5ex plus 1ex minus .2ex}{1.5ex plus 1ex
minus .2ex}{\kai}{\parindent}{\song\bfseries}{}{1em}{}
\numberwithin{equation}{section}

\newcommand{\be}{\begin{equation}}
\newcommand{\ee}{\end{equation}}

\newcommand{\ptl}{\partial}

\newcommand{\eps}{\varepsilon}

\newcommand{\lam}{\lambda}
\numberwithin{equation}{section}

\newtheorem{lemma}{Lemma}[section]

\newtheorem{thm}{Theorem}[section]
\newtheorem{pro}{Proposition}[section]

\allowdisplaybreaks[4]
\begin{document}
\title{{\textbf{Classification of solutions for some mixed order elliptic system}}}
\author{Genggeng Huang\footnote{genggenghuang@fudan.edu.cn} and Yating Niu\footnote{ytniu19@fudan.edu.cn}}
\date{}
\maketitle
\begin{center}
School of Mathematical Sciences, Fudan University, Shanghai, China
\end{center}
\begin{abstract}
In this paper, we classify the solution of the following mixed-order conformally invariant system with coupled nonlinearity in $ \mathbb{R}^4$:
\begin{equation}\label{b74}
\left\{
\begin{aligned}
& -\Delta u(x) = u^{p_1}(x) e^{q_1v(x)}, \quad x\in \mathbb{R}^4,\\
& (-\Delta)^2 v(x) = u^{p_2}(x) e^{q_2v(x)}, \quad x\in \mathbb{R}^4,
\end{aligned}
\right.
\end{equation}
where $ 0\leq p_1 < 1$, $ p_2 >0$, $ q_1 > 0$, $ q_2 \geq 0$, $ u>0$ and satisfies
$$ \int_{\mathbb{R}^4}  u^{p_1}(x) e^{q_1v(x)} dx < \infty,\quad \int_{\mathbb{R}^4}  u^{p_2}(x) e^{q_2 v(x)} dx < \infty.$$
Under additional assumptions \eqref{b20} or \eqref{b26}, we study the asymptotic behavior of the solutions to system \eqref{b74} and we establish the equivalent integral formula for \eqref{b74}. By using the method of moving spheres, we obtain the classification results of the solutions in \eqref{b74}.
\end{abstract}

\section{Introduction}
In this paper, we study the classification results of the following elliptic system
\begin{equation}\label{b1}
\left\{
\begin{aligned}
& -\Delta u(x) = u^{p_1}(x) e^{q_1v(x)}, \quad x\in \mathbb{R}^4,\\
& (-\Delta)^2 v(x) = u^{p_2}(x) e^{q_2v(x)}, \quad x\in \mathbb{R}^4,
\end{aligned}
\right.
\end{equation}
where $ 0\leq p_1 < 1$, $ p_2 >0$, $ q_1 > 0$, $ q_2 \geq 0$, $ u>0$ and $ (u,v)$ satisfies
\be\label{b21}
\ \int_{\mathbb{R}^4}  u^{p_1}(x) e^{q_1v(x)} dx < \infty,  \quad \int_{\mathbb{R}^4}  u^{p_2}(x) e^{q_2 v(x)} dx < \infty \quad \text{and} \quad v(x) = o(|x|^2). \tag{G}
\ee
System \eqref{b1} is conformal invariant for $(p_1,q_1), (p_2,q_2)$ satisfying \eqref{b2}.

\par  For a single equation, there are two simple conformal invariant cases.

\be\label{b10}
-\Delta u= u^{\frac{n+2}{n-2}}, \quad \text{in} \quad \mathbb{R}^n,\quad n\geq 3.
\ee
Or
\be\label{b11}
-\Delta u= e^u, \quad \text{in} \quad \mathbb{R}^2.
\ee
In \cite{GiNiNir79}, Gidas, Ni and Nirenberg studied the symmetry of the solutions of \eqref{b10}.
They proved that all positive solutions satisfying $ u= O(|x|^{2-n})$ at $ \infty$ are radially symmetric about some point. Next, Caffarelli, Gidas and Spruck \cite{CaGiSpr89} removed the growth assumption and proved the similar result. Chen and Li \cite{ChenLi91} also established the classification results for \eqref{b11} and the following equation
under the assumption $ \int_{\mathbb{R}^2} e^{u(x)}dx < \infty$.
\par The natural generalizations of \eqref{b10} and \eqref{b11} are the following higher order equations
\begin{equation}\label{b30}
\left\{
\begin{aligned}
& (-\Delta)^{\frac{n}{2}} u =(n-1)!e^{nu} \quad \text{in} \quad \mathbb{R}^n,\\
& \int_{\mathbb{R}^n} e^{nu(x)}dx < \infty \quad \text{and} \quad u=o(|x|^2) \ \text{at} \ \infty,
\end{aligned}
\right.
\end{equation}
and
\begin{equation}\label{b35}
\left\{
\begin{aligned}
& (-\Delta)^{k} u =u^{\frac{n+2k}{n-2k}} \quad \text{in} \quad \mathbb{R}^n, n>2k,\\
& u>0,\quad \text{in} \quad \mathbb{R}^n.
\end{aligned}
\right.
\end{equation}
Zhu \cite{Zhu} obtained the classification results to \eqref{b30} with $ n=3$. For $ n=4$, Lin \cite{Lin} classified the fourth order equation of \eqref{b30}. Independently, Xu \cite{Xu} proved the same result with $ u=o(|x|^2)$ replaced by
\[
 \lim\limits_{|x| \rightarrow \infty} \Delta u = 0.
\]
Lin \cite{Lin} also considered the classification results of \eqref{b35} with $ k=2$. In the case that $ n$ is an even integer, Wei and Xu \cite{XuWei} classified the smooth solutions of \eqref{b30} and the smooth positive solutions of \eqref{b35}.
The classification of solutions for higher order elliptic equations is much more difficult than second order elliptic equations. For more literatures on higher order equations, please refer to \cite{ChangYang,ChenFangLi,HuangLi,Xu05}.

However, the non-linear terms in generalizations \eqref{b30} and \eqref{b35} are  either purely exponential or purely polynomial. Recently, Yu  \cite{Yu} made an attempt to investigate the generalization of \eqref{b10}, \eqref{b11} to systems with mixed non-linear terms.
In particular, Yu \cite{Yu} consider the following  elliptic system
\begin{equation}\label{b81}
\left\{
\begin{aligned}
&- \Delta u(x) = e^{3 v(x)}, \quad x \in  \mathbb{R}^4, \\
& (- \Delta)^2 v(x) = u^4(x), \quad x \in  \mathbb{R}^4,
\end{aligned}
\right.
\end{equation}
for $u>0$ which corresponds to the case $ p_1 = q_2 =0$, $ q_1 = 3$, $ p_2 = 4$ of the system \eqref{b1}. In \cite{Yu}, Yu classified the solution of the above system \eqref{b81} under the assumptions
\be\label{b20}
\int_{\mathbb{R}^4} e^{3 v(x)} dx < \infty, \quad \int_{\mathbb{R}^4} u^4(x) dx < \infty, \quad \text{and} \quad v(x) = o(|x|^2). \tag{H1}
\ee
Dai et al. \cite{Daiw} obtained the classification results of the following conformally invariant system with mixed order and exponentially increasing nonlinearity in $ \mathbb{R}^2$
\begin{equation}
\left\{
\begin{aligned}
&(- \Delta)^{\frac{1}{2}} u(x) = e^{p v(x)}, \quad x \in  \mathbb{R}^2, \\
& - \Delta v(x) = u^4(x), \quad x \in  \mathbb{R}^2,
\nonumber
\end{aligned}
\right.
\end{equation}
where $ p \in (0,+\infty)$, $ u \geq 0$ and $ \int_{\mathbb{R}^2}  u^{4}(x) dx < \infty$.  Later, Guo and Peng \cite{Guopeng} considered
\begin{equation}\label{gp1}
\left\{
\begin{aligned}
&(- \Delta)^{\frac{1}{2}} u(x) = u^{p_1}(x) e^{q_1 v(x)}, \quad x \in  \mathbb{R}^2, \\
& - \Delta v(x) = u^{p_2}(x)e^{q_2 v(x)}, \quad x \in  \mathbb{R}^2,
\nonumber
\end{aligned}
\right.
\end{equation}
where $ u>0$, $ 0\leq p_1 < \frac{1}{1+K}$, $ p_2 >0$, $ q_1 > 0$, $ q_2 \geq 0$ and $ \int_{\mathbb{R}^2}  u^{p_2}(x)e^{q_2 v(x)} dx < \infty$. Under the assumptions, $ u(x) = O(|x|^K)$ at $ \infty$ for some $ K\geq 1$ arbitrarily large and $ v^+(x) = O(\ln |x|)$ if $ q_2 > 0$ at $ \infty$. It should be pointed out that the coupling of polynomial non-linearity and exponential non-linearity in \eqref{gp1} makes the problem more complicated.
Other classification results for system can be found in \cite{ChenLi05,ChenLi2009,ChenLiOu05,Daiw,GuoLiu08,Guopeng,YangYu} and the references therein.

In this paper, we study the Liouville theorem and classification results for elliptic system \eqref{b1}. Our result is the following

\begin{thm}\label{b59}
Suppose that $ 0\leq p_1 < 1$, $ p_2 >0$, $ q_1 > 0$ and $ q_2 \geq 0$ such that the following conditions hold:
\begin{equation}\label{b2}
\left\{
\begin{aligned}
& \frac{3-p_1}{q_1}= \frac{4-p_2}{q_2}, \quad \text{if} \ q_2>0; \\
& p_2 =4 , \quad \text{if} \ q_2=0.
\end{aligned}
\right.
\end{equation}
Every $ C^4$ solution of \eqref{b1} satisfies \eqref{b20} and \eqref{b21} has the following form
\be\label{b80}
u(x)= \frac{C_1(\eps)}{|x-x_0|^2 + \eps^2}, \quad v(x)= \frac{3-p_1}{q_1}\ln \left( \frac{C_2(\eps)}{|x-x_0|^2 + \eps^2} \right),
\ee
where $ C_1$, $ C_2$ are two positive constants depending only on $ \eps$ and $ x_0$ is a fixed point in $ \mathbb{R}^4$. Assumption \eqref{b20} can be replaced by
\be\label{b26}
 \int_{\mathbb{R}^4}  u^{\tau}(x) dx < \infty, \quad   \tau > p_2, \quad \text{and} \quad v(x)=O(\ln |x|) \quad \text{as} \quad  |x| \rightarrow \infty, \tag{H2}
\ee
where $ 0 \leq p_1 < \min\{\frac{\tau}{3}, 1\}$.
\end{thm}

Theorem \ref{b59} is proved by the moving spheres method, which is a variant of the moving planes method. In recent years, there has been great interest in using the method of moving spheres and moving planes to classify the solutions of equation. It is a very powerful tool to study the symmetry of solutions. For more results, we refer to \cite{ChenLi91,ChenLi2009,ChenLiOu05,ChenLiOu06,Huang,LiyZhang,Niu,Yu13}.
Interested readers may see the book of Chen and Li \cite{ChenLiBook}.

Theorem \ref{b59} generalizes the result in \cite{Yu} and \cite{Guopeng}. Our conditions \eqref{b20} and \eqref{b21} in the case of $ p_1 = q_2 =0$, $ q_1 = 3$, $ p_2 = 4$ are the same as the assumptions in \cite{Yu}. In \cite{Guopeng}, Guo and Peng assumed that $ v^+(x) = O(\ln |x|)$. Inspired by this, we give an alternative assumption \eqref{b26}. We prove that the classification result whether $ u$ satisfies \eqref{b21}\eqref{b20} or \eqref{b21}\eqref{b26}.

This paper is organized as follows. In Section 2, we give the integral representation formula for $ (u,v)$. At the same time, we give the asymptotic behavior of $ v$ at $ \infty$ under two assumptions. In Section 3, we use the moving spheres method to prove that the slow decay of $ v$ can't occur. That is, $ \alpha < \mu$ can't occur. In Section 4, we prove that the fast decay can't occur. In Section 5, we obtain the precise decay of $ v$ and we give the proof of Theorem \ref{b1}.

\section{Preliminaries}
In this section, we establish some estimations on the decay of the solution. We denote
\begin{equation}\label{b3}
\alpha = \frac{1}{4 \omega_3} \int_{\mathbb{R}^4}  u^{p_2}(x) e^{q_2 v(x)} dx,
\end{equation}
where $ \omega_3$ is the volume of the unit sphere in $ \mathbb{R}^4$. We also define
\begin{equation}\label{b4}
k(x)=\frac{1}{4 \omega_3} \int_{\mathbb{R}^4} [\ln|x-y|-\ln(|y|+1)]u^{p_2}(y) e^{q_2 v(y)} dy,
\end{equation}
then it is easy to see that $ k(x)$ is well-defined since $ \int_{\mathbb{R}^4}  u^{p_2}(x) e^{q_2 v(x)} dx < \infty$ and $ u\in C^4(\mathbb{R}^4)$. Moreover, it can be verified that $ k(x)$ satisfies
\begin{equation}\label{b5}
\Delta k(x)= \frac{1}{2 \omega_3} \int_{\mathbb{R}^4} \frac{ u^{p_2}(y) e^{q_2 v(y)} }{|x-y|^2}dy
\end{equation}
and
\begin{equation}\label{b6}
\Delta^2 k(x) = - u^{p_2}(x) e^{q_2 v(x)}.
\end{equation}
Now, we derive some estimations on the decay of $ k(x)$ at infinity.

\begin{lemma}\label{b49}
Let $ k(x)$ be defined as equation \eqref{b4}. Suppose that $ (u,v)$ satisfies \eqref{b21}. Then, we have
\begin{equation}\label{b7}
k(x) \leq \alpha \ln|x|
\end{equation}
for $ |x|$ large enough, where $ \alpha$ is defined by equation \eqref{b3}.
\end{lemma}

\noindent\emph{Proof.} Set
\[
A_1=\left\{y \in \mathbb{R}^4 | |y-x|\leq \frac{|x|}{2} \right\}\ \text{and} \ \ A_2=\left\{y \in \mathbb{R}^4 | |y-x| > \frac{|x|}{2}\right\}.
\]
For $ y\in A_1$, we have $ |x-y|\leq \frac{|x|}{2} \leq |y| < |y|+1$, hence
\[
\ln \frac{|x-y|}{|y|+1} \leq 0,
\]
which further implies
\[
k(x) \leq \frac{1}{4 \omega_3} \int_{A_2} \ln\frac{|x-y|}{|y|+1} u^{p_2}(y) e^{q_2 v(y)} dy.
\]
For $ |x|\geq 2$, we find $ |x-y| \leq |x|+|y| \leq |x|(1+|y|)$ and
\[
\ln \frac{|x-y|}{|y|+1} \leq \ln|x|.
\]
For $ |x|\geq 2$, we conclude that
\begin{align*}
k(x)& \leq \frac{1}{4 \omega_3} \ln|x| \int_{A_2}  u^{p_2}(y) e^{q_2 v(y)} dy \\
& \leq \alpha \ln|x|.
\end{align*}
\hfill $ \Box$

\begin{lemma}\label{b14}
If $ (u,v)$ is a solution of problem \eqref{b1} and it satisfies \eqref{b20} or \eqref{b26}, then
\begin{equation}\label{b8}
u(x) = \frac{1}{2 \omega_3} \int_{\mathbb{R}^4} \frac{u^{p_1}(y) e^{q_1 v(y)}}{|x-y|^2} dy.
\end{equation}
\end{lemma}
\noindent \emph{Proof.} By \eqref{b1}, we get
\[
-\Delta u(x) = u^{p_1}(x) e^{q_1 v(x)}.
\]
By multiplying the equation by $ \frac{1}{|x|^2} - \frac{1}{R^2}$ and integrating over $ B_R$, we obtain
\[
\int_{B_R} \left(\frac{1}{|x|^2} - \frac{1}{R^2} \right) \Delta u(x) dx = - \int_{B_R} \left(\frac{1}{|x|^2} - \frac{1}{R^2} \right) u^{p_1}(x) e^{q_1 v(x)} dx.
\]
By a simple calculation, one gets
\begin{align*}
\int_{B_R \setminus B_{\eps}} \left(\frac{1}{|x|^2} - \frac{1}{R^2} \right) \Delta u(x) dx & = \int_{\ptl(B_R \setminus B_{\eps})}\left(\frac{1}{|x|^2} - \frac{1}{R^2} \right) \nabla u \cdot \nu dS - \int_{B_R \setminus B_{\eps} } \nabla u \cdot  \nabla\left({\frac{1}{|x|^2}}\right) dx \\
& = -\int_{\ptl B_{\eps}} \left(\frac{1}{|x|^2} - \frac{1}{R^2} \right) \nabla u \cdot \nu dS - \int_{\ptl(B_R \setminus B_{\eps})} u \frac{ \ptl \left(\frac{1}{|x|^2}  \right)}{\ptl \nu} dS \\
& =O(\eps)+\frac{2}{R^3} \int_{\ptl B_R} u dS - \frac{2}{\eps^3} \int_{\ptl B_\eps} u dS,
\end{align*}
where $ \nu$ is the outward pointing unit normal vector.  Letting $ \eps \rightarrow 0$, we find
\[
u(0)= \frac{1}{\omega_3 R^3} \int_{\ptl B_R} u dS + \frac{1}{2\omega_3} \int_{B_R} \left(\frac{1}{|x|^2} - \frac{1}{R^2} \right) u^{p_1}(x) e^{q_1 v(x)} dx.
\]
That is
\[
u(x)= \frac{1}{\omega_3 R^3} \int_{\ptl B_R(x)} u(y) dS + \frac{1}{2\omega_3} \int_{B_R(x)} \left(\frac{1}{|y-x|^2} - \frac{1}{R^2} \right) u^{p_1}(y) e^{q_1 v(y)} dy.
\]
Since $ u(x) > 0$, we have
\[
u(x) \geq \frac{1}{2\omega_3} \int_{B_R(x)} \left(\frac{1}{|y-x|^2} - \frac{1}{R^2} \right) u^{p_1}(y) e^{q_1 v(y)} dy.
\]
Letting $ R \rightarrow \infty$, we conclude that
\begin{equation*}
u(x) \ge \frac{1}{2 \omega_3} \int_{\mathbb{R}^4} \frac{u^{p_1}(y) e^{q_1 v(y)}}{|x-y|^2} dy.
\end{equation*}
Since $\tilde u(x)= u(x)- \frac{1}{2 \omega_3} \int_{\mathbb{R}^4} \frac{u^{p_1}(y) e^{q_1 v(y)}}{|x-y|^2} dy$ is harmonic, by Liouville theorem for non-negative harmonic functions, $\tilde u\equiv C_0\ge 0$.


 By our assumption \eqref{b20} or \eqref{b26}, we obtain $ C_0 = 0$. \hfill $ \Box$

\begin{lemma}\label{b13}
Suppose that $ (u,v)$ is a solution of problem \eqref{b1} with \eqref{b21} \eqref{b20} or \eqref{b21} \eqref{b26}. Then we have
\begin{equation}\label{b9}
\Delta v(x) = - \frac{1}{2 \omega_3} \int_{\mathbb{R}^4} \frac{u^{p_2}(y) e^{q_2 v(y)}}{|x-y|^2} dy -\tilde{C}
\end{equation}
for some $ \tilde{C}\geq 0$.
\end{lemma}
\noindent\emph{Proof.} Let $ k(x)$ be as \eqref{b4} and $ h(x)=k(x) + v(x)$. Then $ \Delta^2 h = 0$, i.e., $ \Delta h$ is harmonic in $ \mathbb{R}^4$. For any $ x_0 \in \mathbb{R}^4$, we infer from the mean value theorem of harmonic function that
\begin{align*}
\Delta h(x_0)& = \frac{4}{\omega_3 s^4} \int_{B_s(x_0)} \Delta h(x) dx \\
& = \frac{4}{\omega_3 s^4} \int_{\ptl B_s(x_0)} \frac{\ptl h}{\ptl \nu} dS.
\end{align*}
Integrating from $ 0$ to $ r$, we get
\[
\frac{r^2}{8} \Delta h(x_0) = \frac{1}{\omega_3 r^3} \int_{\ptl B_r(x_0)} h(y) dS - h(x_0),
\]
which further implies
\begin{align*}
e^{\frac{q_1 r^2}{8} \Delta h(x_0)} & = e^{-q_1 h(x_0)}e^{q_1 -\hspace{-0.5em}\int_{\ptl B_r(x_0)} h(y) dS} \\
& \leq  e^{-q_1 h(x_0)} -\hspace{-1.1em}\int_{\ptl B_r(x_0)} e^{q_1 h(y)}dS \\
& = e^{-q_1 h(x_0)} \frac{1}{\omega_3 r^3} \int_{\ptl B_r(x_0)} e^{q_1 h(y)}dS.
\end{align*}
From the integral representation formula \eqref{b8} for $ u$, one gets
\begin{equation}\label{b12}
\begin{aligned}
u(x) &\geq \frac{1}{2 \omega_3} \int_{ |y| < \frac{|x|}{2}} \frac{u^{p_1}(y) e^{q_1 v(y)}}{|x-y|^2} dy \\
& \geq \frac{2}{9 \omega_3 |x|^2} \int_{ |y| < 1} u^{p_1}(y) e^{q_1 v(y)} dy =: \frac{C}{|x|^2},
\end{aligned}
\end{equation}
for any $ |x|\ge 2$.  Since $ k(x) \leq \alpha \ln|x|$, then
\[
e^{q_1 h(x)} \leq e^{q_1 v(x)}e^{q_1 \alpha \ln|x|} = |x|^{\alpha q_1} e^{q_1 v(x)}.
\]
By the assumption $ \int_{\mathbb{R}^4}  u^{p_1}(x) e^{q_1 v(x)} dx < \infty$, we have $ |x|^{- \alpha q_1}e^{q_1 h(x)} u^{p_1}(x)$ is integrable in $ \mathbb{R}^4 \setminus B_2(0)$. Since \eqref{b12}, we get $ r^{- (\alpha q_1 + 2p_1)}\int_{\ptl B_r(x_0)} e^{q_1 h(x)} dS$ is integrable in $ [R,+\infty)$, where $ R> |x_0|+2$. Therefore,
\[
e^{\frac{q_1 r^2}{8} \Delta h(x_0)} r^{3 - (\alpha q_1 + 2p_1)} \leq C r^{- (\alpha q_1 + 2p_1)} \int_{\ptl B_r(x_0)} e^{q_1 h(y)}dS
\]
is integrable in $ [R,+\infty)$. In particular, we obtain $ \Delta h(x_0) \leq 0$, $ \forall x_0 \in \mathbb{R}^4$. Hence, we infer from the Liouville theorem for harmonic function that $ \Delta h(x) = -\tilde{C}$ for some $ \tilde{C}\geq 0$. We conclude that
\[
\Delta v(x) = - \frac{1}{2 \omega_3} \int_{\mathbb{R}^4} \frac{u^{p_2}(y) e^{q_2 v(y)}}{|x-y|^2} dy -\tilde{C},
\]
for some $ \tilde{C}\geq 0$. \hfill $ \Box$

By \eqref{b20} or \eqref{b26}, we find $ v(x)=o(|x|^2)$ at $ \infty$. Then, we obtain the following lemma.

\begin{lemma}\label{b29}
Suppose that $ v(x)=o(|x|^2)$ as $ |x|\rightarrow \infty$, then we have
\be\label{b18}
\Delta v(x) = - \frac{1}{2 \omega_3} \int_{\mathbb{R}^4} \frac{u^{p_2}(y) e^{q_2 v(y)}}{|x-y|^2} dy.
\ee
\end{lemma}
\noindent \emph{Proof.}
By Lemma \ref{b13}, we know that $ \Delta v \le -\tilde C$. If $ \tilde{C} > 0$, let $ w= -v -\frac{\tilde{C}}{8} |x|^2$. Then $\Delta w\ge 0$ and hence $w(0)\le \sup_{\partial B_r} w$, i.e.,
\[
\inf_{\ptl B_r} v \leq v(0) - \frac{\tilde{C}}{8} r^2.
\]
This contradicts the assumption $ v(x)=o(|x|^2)$. This completes the proof of this lemma.  \hfill $ \Box$

\begin{lemma}\label{b22}
There exists a constant $ C >0$, such that
\[
u(x) \leq C, \quad v(x) \leq C \quad \text{and} \quad | \Delta v(x)| \leq C \quad \text{in} \quad \mathbb{R}^4,
\]
under the assumptions \eqref{b21} and \eqref{b20}.
\end{lemma}
Before we prove the above Lemma \ref{b22}, we need to recall a Brezis-Merle type result, for its proof, see \cite{Lin}.

\begin{lemma}\label{b23}
Suppose that $ h(x)$ satisfies
\begin{equation}
\left\{
\begin{aligned}
& \Delta^2 h(x) =g(x) \quad \text{in} \ B_R \subset \mathbb{R}^4, \\
& h(x)=\Delta h(x) = 0 \ \text{on} \ \ptl B_R,
\nonumber
\end{aligned}
\right.
\end{equation}
with $ g\in L^1(B_R)$, then for any $ \delta \in (0,32\pi^2)$, there exists a $ C_{\delta} >0$, such that
\[
\int_{B_R} e^{\frac{\delta |h|}{\|g\|_{L^1}}} \leq C_{\delta}R^4.
\]
\end{lemma}
\noindent\emph{Proof of Lemma \ref{b22}.} Let $ h(y)$ be defined as
\begin{equation}
\left\{
\begin{aligned}
& \Delta^2 h(y) = u^{p_2}(y) e^{q_2 v(y)} \quad \text{in} \ B_4(x), \\
& h(y) = \Delta h(y) = 0 \ \text{on} \ \ptl B_4(x).
\nonumber
\end{aligned}
\right.
\end{equation}
Since $ \int_{\mathbb{R}^4}  u^{p_2}(y) e^{q_2 v(y)} dy < \infty$, then it follows from Lemma \ref{b23} that there exists an $ L>0$ such that $ \int_{B_4(x)} e^{\theta |h|} dy \leq C$,  $ \forall x\in B^c_L$, where $ \theta = \max\{12 q_1 , 3q_2\}$ and $ C$ is a constant independent of $ x$.

Set $ q(y)=v(y)-h(y)$ in $ B_4(x)$, then $ q(y)$ satisfies
\begin{equation}
\left\{
\begin{aligned}
& \Delta^2 q(y) = 0 \quad \text{in} \ B_4(x), \\
& q(y) = v(y), \ \Delta q(y) = \Delta v(y) \ \text{on} \ \ptl B_4(x).
\nonumber
\end{aligned}
\right.
\end{equation}
Next, let $ \tilde{q}(y) = -\Delta q(y)$, then $ \tilde{q}(y)$ is harmonic and non-negative in $ B_4(x)$. Therefore, we infer from Harnack inequality and mean value theorem that
\[
\tilde{q}(y) \leq C_2 \tilde{q}(x) = C_2 -\hspace{-1.1em}\int_{\ptl B_4(x)} \tilde{q}(z)dS = C_2 -\hspace{-1.1em}\int_{\ptl B_4(x)} - \Delta v(z) dS, \qquad \forall y\in B_2(x),
\]
where the constant $ C_2$ is independent of $ x$ and $ y$. On the other hand, we note that $ v(x)$ satisfies
\[
 (-\Delta)^2 v(x) = u^{p_2}(x) e^{q_2 v(x)} \quad \text{in} \ \mathbb{R}^4.
\]
Similarly to the calculation of Lemma \ref{b14}, one gets
\begin{align*}
 - -\hspace{-1.1em}\int_{\ptl B_r(x)} \Delta v(y) dS & = - \Delta v(x) - \frac{1}{2 \omega_3} \int_{B_r(x)}\left(\frac{1}{|x-y|^2}-\frac{1}{r^2}  \right) u^{p_2}(y) e^{q_2 v(y)} dy \\
& = \frac{1}{2 \omega_3} \int_{\mathbb{R}^4 \setminus B_r(x)} \frac{u^{p_2}(y) e^{q_2 v(y)}}{|x-y|^2}dy + \frac{1}{2 \omega_3 r^2} \int_{ B_r(x)} u^{p_2}(y) e^{q_2 v(y)}dy.
\end{align*}
Taking $ r=4$, we get
\[
- -\hspace{-1.1em}\int_{\ptl B_4(x)} \Delta v(y) dS < C.
\]
In particular, we have $ \tilde{q}(y) \leq C$ for $ y \in B_2(x)$.

On the other hand, since $ q(y)$ satisfies
\begin{equation}
\left\{
\begin{aligned}
& \Delta q(y) = -\tilde{q}(y) \quad \text{in} \ B_4(x), \\
& q(y) = v(y) \quad \text{on} \ \ptl B_4(x),
\nonumber
\end{aligned}
\right.
\end{equation}
then
\be\label{b25}
\sup_{B_1(x)} q \leq C\{\|q^+\|_{L^2(B_2(x))}+ \|\tilde{q}\|_{L^\infty (B_2(x))}\}.
\ee
Since $ q = v - h$, then
\[
q^+ \leq v^+ + |h|,
\]
which further implies
\be\label{b24}
\int_{B_2(x)} (q^+)^2 dy \leq 2 \int_{B_2(x)} e^{2 v^+} + |h|^2 dy \leq C \int_{B_2(x)} (1+ e^{3v})dy + C \int_{B_2(x)} e^{3|h|} dy \leq C.
\ee
Substituting \eqref{b24} into \eqref{b25}, we deduce that $ \sup\limits_{B_1(x)} q$ is bounded. Therefore, we find $ v(y) \leq C+ |h(y)|$ for $ y \in B_1(x)$, which further implies
\[
\int_{ B_1(x)} e^{\theta v(y)} dy \leq C \int_{ B_1(x)} e^{\theta |h(y)|} dy \leq C.
\]
By Lemma \ref{b14}, one gets
\begin{align*}
u(x) & = \frac{1}{2 \omega_3} \int_{\mathbb{R}^4} \frac{u^{p_1}(y) e^{q_1 v(y)}}{|x-y|^2} dy \\
& = \frac{1}{2 \omega_3} \int_{B_1^c (x)} \frac{u^{p_1}(y) e^{q_1 v(y)}}{|x-y|^2} dy + \frac{1}{2 \omega_3} \int_{B_1(x)} \frac{u^{p_1}(y) e^{q_1 v(y)}}{|x-y|^2} dy \\
& \leq C + C\left(\int_{B_1(x)} \frac{1}{|x-y|^3} dy \right)^{\frac{2}{3}} \left(\int_{B_1(x)} u^4(y) dy \right)^{\frac{p_1}{4}} \left(\int_{B_1(x)} e^{s q_1 v(y)} dy \right)^{\frac{1}{s}} \\
& \leq C,
\end{align*}
where $ \frac{2}{3} + \frac{p_1}{4} + \frac{1}{s} = 1$, $ 0 \leq p_1 < 1$, $ 3 \leq s <12$. Hence, we obtain
\begin{align*}
- \Delta v(x) & = \frac{1}{2 \omega_3} \int_{\mathbb{R}^4} \frac{u^{p_2}(y) e^{q_2 v(y)}}{|x-y|^2} dy \\
& = \frac{1}{2 \omega_3} \int_{B_1^c (x)} \frac{u^{p_2}(y) e^{q_2 v(y)}}{|x-y|^2} dy + \frac{1}{2 \omega_3} \int_{B_1(x)} \frac{u^{p_2}(y) e^{q_2 v(y)}}{|x-y|^2} dy \\
& \leq C +  C \left(\int_{B_1(x)} \frac{1}{|x-y|^3} dy \right)^{\frac{2}{3}} \left(\int_{B_1(x)} e^{3 q_2 v(y)} dy \right)^{\frac{1}{3}} \\
& \leq C.
\end{align*}

Since $ - \Delta v > 0$, we get $ |\Delta v| \leq C$, where $ C > 0$ is a constant. Let $ f = - \Delta v$. Recalling that $ \Delta v \in C^2$, we take $ v_1$ to be a solution of
\begin{equation}
\left\{
\begin{aligned}
- \Delta v_1 = f \quad \text{in}  \ B_1(x) \\
v_1 = 0  \quad \text{on}  \ \ptl B_1(x).
\nonumber
\end{aligned}
\right.
\end{equation}

From elliptic theory we know that $ |v_1| \leq C_1$, for some constant $ C_1$ depending on $ C$ but not depending on $ x$. Let $ v_2 = v - v_1$. Then $ \Delta v_2 = 0$ in $ B_1(x)$. By the mean value theorem of harmonic functions, we observe that
\begin{align*}
v_2^+ (x)& \leq C_2 \int_{B_1(x)} v_2^+ dx \\
& \leq C_2 \left( \int_{B_1(x)}(1+ e^{3v})dx + \int_{B_1(x)} |v_1| dx \right) \leq C_2.
\end{align*}
Therefore, we conclude that $ v \leq C$. \hfill $ \Box$

\begin{lemma}\label{b77}
Suppose that $ (u,v)$ satisfies \eqref{b21} and \eqref{b20}. Then  $\forall \eps >0$, there exists an $R_\varepsilon>0$ such that
\[
k(x) \geq (\alpha - \eps) \ln |x|,\quad \forall |x|\ge R_\varepsilon.
\]
\end{lemma}
\noindent \emph{Proof.} Let $ A_1=\{y \in \mathbb{R}^4||y|\leq R_0\}$. Then we can choose $ R_0$ large enough, such that
\begin{equation}\label{b47}
\frac{1}{4 \omega_3} \int_{A_1} [\ln|x-y|-\ln(|y|+1)]u^{p_2}(y) e^{q_2 v(y)} dy \geq \left(\alpha - \frac{\eps}{2}\right) \ln|x|.
\end{equation}
Let $ A_2=\left\{y \in \mathbb{R}^4 ||y-x|\leq \frac{|x|}{2}, |y|\geq R_0 \right\}$ and $ A_3=\left\{y \in \mathbb{R}^4||y-x|> \frac{|x|}{2}, |y|\geq R_0 \right\}$, then
\begin{equation}\label{b72}
\begin{aligned}
& \int_{A_2} [\ln|x-y|-\ln(|y|+1)]u^{p_2}(y) e^{q_2 v(y)} dy \\
& \geq \int_{B_1(x)}\ln|x-y|u^{p_2}(y) e^{q_2 v(y)} dy - \int_{A_2} \ln(|y|+1)u^{p_2}(y) e^{q_2 v(y)} dy \\
& \geq -C - \ln(2|x|)\int_{A_2} u^{p_2}(y) e^{q_2 v(y)} dy \\
& \geq -C - \frac{\eps}{4} \ln|x|.
\end{aligned}
\end{equation}
Finally, we estimate
\[
\int_{A_3}[\ln|x-y|-\ln(|y|+1)]u^{p_2}(y) e^{q_2 v(y)} dy.
\]
If $ y\in A_3$ and $ |y|\leq 2|x|$, then we have $ |x-y| > \frac{|x|}{2} \geq \frac{|y|}{4}$. If $ y\in A_3$ and $ |y| > 2|x|$, then we find $ |x-y| \geq |y|-|x| \geq |y| - \frac{|y|}{2} = \frac{|y|}{2}$. That is, in both cases, we observe that
\[
\frac{|x-y|}{|y|} \geq \frac{1}{4}
\]
or
\[
\frac{|x-y|}{|y|+1} \geq \frac{1}{8}.
\]
Hence, it is clear that
\begin{equation}\label{b73}
\int_{A_3}[\ln|x-y|-\ln(|y|+1)]u^{p_2}(y) e^{q_2 v(y)} dy \geq \ln \frac{1}{8} \int_{A_3} u^{p_2}(y) e^{q_2 v(y)} dy \geq - \frac{\eps}{4} \ln|x|
\end{equation}
for $ |x|$ enough large.
Finally, we infer from \eqref{b47}-\eqref{b73} that
\[
k(x)\geq (\alpha -\eps)\ln|x|.
\]
\hfill $ \Box$

\begin{lemma}\label{b50}
If $ (u,v)$ is a solution of problem \eqref{b1}, then we have
\begin{equation}
\lim\limits_{|x|\rightarrow \infty} \frac{k(x)}{\ln|x|} = \alpha
\end{equation}
under the assumptions \eqref{b21} and \eqref{b26}.
\end{lemma}
\noindent \emph{Proof.} We need to show that
\[
\lim_{|x|\rightarrow \infty} \int_{\mathbb{R}^4} \frac{\ln |x-y| - \ln (|y|+1) - \ln |x|}{\ln |x|} u^{p_2}(y) e^{q_2 v(y)} dy = 0.
\]
Consider $ |x|$ large enough, and set
\[
A_1 = \{y \in \mathbb{R}^4||y-x| \leq 1\}, \ A_2 = \{y\in \mathbb{R}^4 ||y-x|>1, |y| \geq \ln |x|\}
\]
and
\[
 A_3 = \{y \in \mathbb{R}^4||y-x|>1, |y| < \ln |x|\}.
\]
For $ y\in A_1$, one gets
\begin{align*}
& \left| \int_{B_1(x)} \frac{\ln |x-y| - \ln (|y|+1) - \ln |x|}{\ln |x|} u^{p_2}(y) e^{q_2 v(y)} dy \right| \\
& \leq \int_{B_1(x)} \frac{\ln (\frac{1}{|x-y|}) + \ln (|x|+2) + \ln |x|}{\ln |x|} u^{p_2}(y) e^{q_2 v(y)} dy \\
& \leq 3 \int_{B_1(x)} u^{p_2}(y) e^{q_2 v(y)} dy + \frac{J_1}{\ln |x|}.
\end{align*}
Let
\[
A_1^{'} = \{ y| |y-x| \leq |x|^{-\sigma}\}, \quad A_1^{''} = \{y \in A_1| |y-x| > |x|^{-\sigma} \}.
\]
Next, we give the upper bound of $ J_1$
\begin{align*}
J_1 & =  \int_{B_1(x)} \ln \left(\frac{1}{|x-y|}\right) u^{p_2}(y) e^{q_2 v(y)} dy \\
& = \int_{A_1^{'}} \ln \left(\frac{1}{|x-y|}\right) u^{p_2}(y) e^{q_2 v(y)} dy + \int_{A_1^{''}} \ln \left(\frac{1}{|x-y|}\right) u^{p_2}(y) e^{q_2 v(y)} dy \\
& \leq \max_{B_{|x|^{-\sigma}}(x)}(e^{q_2 v(y)}) \left(\int_{A_1^{'}} \left( \ln \left(\frac{1}{|x-y|}\right) \right)^{s} dy \right)^{\frac{1}{s}}  \left(\int_{A_1^{'}} u^{tp_2} dy \right)^{\frac{1}{t}} + \sigma \ln |x| \int_{A_1^{''}} u^{p_2}(y) e^{q_2 v(y)} dy \\
& \leq C\frac{|x|^{C q_2 }}{|x|^{\frac{3\sigma}{s}}}  + \sigma \ln |x| \int_{B_1(x)} u^{p_2}(y) e^{q_2 v(y)} dy,
\end{align*}
where $ \frac{1}{s} + \frac{1}{t} = 1$ and $ t= \frac{\tau}{p_2}$. In getting last inequality of above, we use \eqref{b26} and assume $ v(x) \le C\ln|x|$, for $ |x|$ large. Fix a large $ \sigma$ such that $ \frac{3\sigma}{s} \geq C q_2$.  Then, one has
\begin{align*}
& \left| \int_{B_1(x)} \frac{\ln |x-y| - \ln (|y|+1) - \ln |x|}{\ln |x|} u^{p_2}(y) e^{q_2 v(y)} dy \right| \\
& \leq  C \int_{B_1(x)} u^{p_2}(y) e^{q_2 v(y)} dy  + \frac{C }{\ln|x|} \\
& = o_{|x|}(1).
\end{align*}
For $ y \in A_2$, we find
\[
\frac{1}{4|x|^2} \leq \frac{|x-y|}{|x|(|y|+1)} \leq \frac{1}{|x|} + \frac{1}{|y|+1} < 1
\]
for $ |x|$ large enough. Then
\begin{align*}
& \left| \int_{A_2} \frac{\ln |x-y| - \ln (|y|+1) - \ln |x|}{\ln |x|} u^{p_2}(y) e^{q_2 v(y)} dy \right| \\
& \leq \sup_{y\in A_2} \frac{| \ln |x-y| - \ln (|y|+1) - \ln |x| |}{\ln |x|} \int_{|y| \geq \ln |x|} u^{p_2}(y) e^{q_2 v(y)} dy \\
& \leq \frac{2 \ln 2 + 2 \ln |x|}{\ln|x|} \int_{|y| \geq \ln |x|} u^{p_2}(y) e^{q_2 v(y)} dy \\
& = o_{|x|}(1).
\end{align*}
For $ y\in A_3$,  one gets
\begin{align*}
& \left| \int_{A_3} \frac{\ln |x-y| - \ln (|y|+1) - \ln |x|}{\ln |x|} u^{p_2}(y) e^{q_2 v(y)} dy \right| \\
& \leq \frac{\max\limits_{y\in A_3}\left| \ln (\frac{|x-y|}{|x|})\right|}{\ln |x|} \int_{A_3} u^{p_2}(y) e^{q_2 v(y)} dy + \frac{1}{\ln |x|} \int_{A_3} \ln (|y|+1) u^{p_2}(y) e^{q_2 v(y)} dy \\
& \leq \frac{C \ln 2}{\ln |x|}  + \frac{C \ln (\ln|x|+1)}{\ln |x|} \\
& = o_{|x|}(1),
\end{align*}
where we have used
\[
\frac 12\le 1-\frac{\ln|x|}{|x|}\le \frac{|x-y|}{|x|} \leq 1+ \frac{\ln |x|}{|x|} \leq \frac 32, \quad \forall y\in A_3.
\]
This completes the proof of this lemma. \hfill $ \Box$

\begin{lemma}\label{b78}
Let $ (u,v)$ be a solution of problem \eqref{b1} with $ v(x) = o(|x|^2)$ as $ |x| \rightarrow \infty$. Then
\begin{equation}\label{b27}
v(x) = - \frac{1}{4 \omega_3} \int_{\mathbb{R}^4} [\ln|x-y|-\ln(|y|+1)]u^{p_2}(y) e^{q_2 v(y)} dy + c_0,
\end{equation}
\begin{equation}\label{b28}
\nabla v(x) =  - \frac{1}{4 \omega_3} \int_{\mathbb{R}^4} \frac{x-y}{|x-y|^2} u^{p_2}(y) e^{q_2 v(y)} dy
\end{equation}
for some constant $c_0\in \mathbb R$.
\end{lemma}
\noindent \emph{Proof.}
Since
\[
\Delta (k(x) + v(x)) = 0
\]
and
\[
k(x) + v(x) = o(|x|^2),
\]
then it follows from Liouville theorem that
\[
k(x) + v(x) = \sum^4_{i=1} c_i x_i +c_0.
\]
Hence we observe that
\[
u^{p_1}(x) e^{q_1 v(x)} = u^{p_1}(x) e^{- q_1 k(x)}e^{q_1 c_0} e^{q_1 \sum^4_{i=1} c_i x_i} \geq C e^{q_1 c_0} |x|^{-(2p_1+ \alpha q_1 )} e^{q_1 \sum^4_{i=1} c_i x_i},
\]
for $ |x|$ enough large. We infer from $ \int_{\mathbb{R}^4}  u^{p_1}(x) e^{q_1v(x)} dx < \infty$ that $ c_1=c_2=c_3=c_4 =0$, which further implies
\[
v(x) =-k(x) +c_0.
\]
This proves equation \eqref{b27}. Equation \eqref{b28} follows from equation \eqref{b27}. \hfill $ \Box$

By Lemma \ref{b49}, Lemma \ref{b77} and Lemma \ref{b50}, we obtain $ k(x) \sim \alpha \ln|x|$ at $ \infty$ whether $ u$ satisfies \eqref{b21}\eqref{b20} or \eqref{b21}\eqref{b26}. By Lemma \ref{b78}, we know that $ v(x) \sim - \alpha \ln|x|$ at $ \infty$, that is,
\begin{equation}\label{1109}
\lim_{|x| \rightarrow \infty} \frac{v(x)}{\ln |x|} = - \alpha.
\end{equation}

\begin{lemma}\label{b15}
Suppose that $ (u,v)$ satisfies \eqref{b21} and \eqref{b26}. Then
\[
u(x) \leq C \quad \text{and} \quad v(x) \leq C \quad \text{in} \ \mathbb{R}^4.
\]
\end{lemma}
\noindent \emph{Proof.} By Lemma \ref{b50} and equation \eqref{b27}, we get
\[
 \lim_{|x| \rightarrow \infty} v(x) = - \infty
\]
and $ v(x) \leq C$. Then
\begin{align*}
u(x) &= \frac{1}{2 \omega_3} \int_{\mathbb{R}^4} \frac{u^{p_1}(y) e^{q_1 v(y)}}{|x-y|^2} dy\\
& = \frac{1}{2 \omega_3} \int_{B_1^c(x)} \frac{u^{p_1}(y) e^{q_1 v(y)}}{|x-y|^2} dy + \frac{1}{2 \omega_3} \int_{B_1(x)} \frac{u^{p_1}(y) e^{q_1 v(y)}}{|x-y|^2} dy \\
& \leq C + C \left(\int_{B_1(x)} \frac{1}{|x-y|^3} dy \right)^{\frac{2}{3}} \left(\int_{B_1(x)} u^{3 p_1} dy \right)^{\frac{1}{3}} \\
& \leq C.
\end{align*}
\hfill $ \Box$

By Lemma \ref{b22} and Lemma \ref{b15}, we get $ u(x) \leq C $ and $  v(x) \leq C$ whenever $(u,v)$ satisfies  \eqref{b21}\eqref{b20} or \eqref{b21}\eqref{b26}.

\section{Slow decay is impossible}
Let
\[
\mu = \frac{6-2p_1}{q_1} \quad \text{and} \quad 2p_2 + \mu q_2 =8.
\]
Let $ \alpha$ be defined as \eqref{b3}, and we prove that $ \alpha = \mu$. Then, we obtain the precise decay of $ v(x)$. In this section, we show that $ \alpha < \mu  $ can't occur. We use the moving spheres method to prove our results.

Let
\begin{equation}\label{b34}
\left\{
\begin{aligned}
& w(x) = \frac{1}{|x|^2} u \left(\frac{x}{|x|^2}\right), \\
& z(x) = v \left(\frac{x}{|x|^2}\right) - \mu \ln|x|.
\end{aligned}
\right.
\end{equation}
be the Kelvin transformation of $ (u,v)$. Then a direct calculation shows that they satisfy the following equation
\begin{equation}
\left\{
\begin{aligned}
& -\Delta w(x) = w^{p_1}(x)e^{q_1 z(x)} , \\
& (-\Delta)^2 z(x) = w^{p_2}(x)e^{q_2 z(x)},
\end{aligned}\quad \text{in}\quad \mathbb R^4\backslash\{0\}.
\nonumber
\right.
\end{equation}
Moreover, one has
\be\label{b33}
-\Delta z(x)=4\frac{x}{|x|^4}\nabla v \left(\frac{x}{|x|^2}\right) - \frac{1}{|x|^4}\Delta v \left(\frac{x}{|x|^2} \right) + \frac{2\mu}{|x|^2}, \qquad x\in \mathbb R^4\backslash\{0\}.
\ee

\begin{lemma} We have
\[
-\Delta z(x)= \frac{1}{2\omega_3} \int_{\mathbb{R}^4} \frac{w^{p_2}(y) e^{q_2 z(y)}}{|x-y|^2}  dy + \frac{2(\mu - \alpha)}{|x|^2}
\]
for $ x\in \mathbb{R}^4 \setminus \{0\}$.
\end{lemma}

\noindent \emph{Proof.} By \eqref{b18} and \eqref{b28}, we know that
\be\label{b31}
\frac{\ptl v}{\ptl x_i}= - \frac{1}{4 \omega_3} \int_{\mathbb{R}^4} \frac{x_i-y_i}{|x-y|^2} u^{p_2}(y) e^{q_2 v(y)}dy
\ee
and
\be\label{b32}
\Delta v(x) = - \frac{1}{2 \omega_3} \int_{\mathbb{R}^4} \frac{u^{p_2}(y) e^{q_2 v(y)}}{|x-y|^2} dy.
\ee
Substituting \eqref{b31} and \eqref{b32} into \eqref{b33}, then
\begin{align*}
-\Delta z(x) & = -\frac{1}{\omega_3} \int_{\mathbb{R}^4} \frac{\sum^4_{i=1}\left(\frac{x_i^2}{|x|^4}-\frac{x_i y_i}{|x|^2 |y|^2}\right)}{|x-y|^2 |y|^6} u^{p_2}\left(\frac{y}{|y|^2}\right) e^{q_2 v\left(\frac{y}{|y|^2}\right) } dy \\
&+\frac{1}{2 \omega_3 |x|^4}\int_{\mathbb{R}^4} \frac{u^{p_2}\left(\frac{y}{|y|^2}\right) e^{q_2 v\left(\frac{y}{|y|^2}\right)}}{\left|\frac{x}{|x|^2}-\frac{y}{|y|^2}\right|^2 |y|^8}dy + \frac{2\mu}{|x|^2} \\
& = -\frac{1}{\omega_3} \int_{\mathbb{R}^4} \frac{|y|^2 - \langle x,y \rangle}{|x-y|^2 |x|^2} w^{p_2}(y) e^{q_2 z(y)} dy + \frac{1}{2 \omega_3} \int_{\mathbb{R}^4} \frac{|y|^2 }{|x-y|^2 |x|^2} w^{p_2}(y) e^{q_2 z(y)} dy + \frac{2\mu}{|x|^2} \\
& = -\frac{1}{2 \omega_3} \int_{\mathbb{R}^4} \frac{|y|^2 - 2\langle x,y \rangle}{|x-y|^2 |x|^2}w^{p_2}(y) e^{q_2 z(y)} dy  + \frac{2\mu}{|x|^2} \\
& = -\frac{1}{2 \omega_3} \int_{\mathbb{R}^4} \frac{|x-y|^2 - |x|^2}{|x-y|^2 |x|^2}w^{p_2}(y) e^{q_2 z(y)} dy  + \frac{2\mu}{|x|^2} \\
& = -\frac{1}{2 \omega_3} \int_{\mathbb{R}^4} \frac{w^{p_2}(y) e^{q_2 z(y)}}{|x|^2} dy + \frac{1}{2 \omega_3}\int_{\mathbb{R}^4} \frac{w^{p_2}(y) e^{q_2 z(y)}}{|x-y|^2} dy + \frac{2\mu}{|x|^2} \\
& = \frac{1}{2 \omega_3} \int_{\mathbb{R}^4} \frac{w^{p_2}(y) e^{q_2 z(y)}}{|x-y|^2} dy + \frac{2(\mu - \alpha)}{|x|^2}
\end{align*}
for $ x\in \mathbb{R}^4 \setminus \{0\}$. \hfill $ \Box$

\begin{lemma}\label{b39} We have
\[
\pm 4 x\cdot \nabla v(x) \geq |x|^2 \Delta v(x) - 2\alpha.
\]
\end{lemma}

\noindent \emph{Proof.} By \eqref{b18} and \eqref{b28}, we find
\begin{align*}
4|x\cdot \nabla v(x)| & \leq \frac{|x|}{\omega_3} \int_{\mathbb{R}^4}\frac{u^{p_2}(y) e^{q_2 v(y)} }{|x-y|}dy \\
& \leq \frac{|x|}{\omega_3}  \left( \int_{\mathbb{R}^4} \frac{u^{p_2}(y) e^{q_2 v(y)} }{|x-y|^2}dy \right)^{\frac{1}{2}} \left( \int_{\mathbb{R}^4} u^{p_2}(y) e^{q_2 v(y)} dy \right)^{\frac{1}{2}} \\
& \leq |x|^2(-\Delta v(x)) + 2\alpha.
\end{align*}
\hfill $ \Box$

Let $ w(x)$ and $ z(x)$ be defined as \eqref{b34}, we define
\[
w_{\lam}(x) = \frac{\lam^2}{|x|^2} w\left(\frac{\lam^2 x}{|x|^2}\right) = \frac{1}{\lam^2} u\left(\frac{x}{\lam^2}\right)
\]
and
\[
z_{\lam}(x) = z\left(\frac{\lam^2 x}{|x|^2}\right) - \mu \ln \frac{|x|}{\lam} = v\left( \frac{x}{\lam^2}\right) - \mu \ln \lam.
\]
Then a direct calculation shows that they satisfy the following equation
\[
-\Delta w_{\lam}(x) = w^{p_1}_{\lam}(x) e^{q_1 z_{\lam}(x)},
\]
\[
-\Delta z_{\lam}(x) = \frac{1}{2 \omega_3} \int_{\mathbb{R}^4} \frac{w^{p_2}_{\lam}(y) e^{q_2 z_{\lam}(y)}}{|x-y|^2} dy
\]
and
\[
\Delta^2 z_{\lam}(x) = w^{p_2}_{\lam}(x) e^{q_2 z_{\lam}(x)}
\]
for any $ x \in \mathbb{R}^4$. Set
\[
S_\lam^w(x) = w(x) - w_{\lam}(x)
\]
and
\[
S_\lam^z(x) = z(x) - z_{\lam}(x).
\]
Then $ S_\lam^w(x)$ satisfies
\[
-\Delta S_\lam^w(x)=w^{p_1}(x) e^{q_1 z(x)} - w^{p_1}_{\lam}(x) e^{q_1 z_{\lam}(x)},\quad \text{in}\quad \mathbb{R}^4 \backslash\{0\}.
\]
Also in $ \mathbb{R}^4\backslash\{0\}$, $ S_\lam^z(x)$ satisfies
\begin{equation}\label{b42}
\begin{aligned}
-\Delta S_\lam^z(x) & = \frac{1}{2 \omega_3} \int_{\mathbb{R}^4} \frac{w^{p_2}(y) e^{q_2 z(y)}}{|x-y|^2} dy - \frac{1}{2 \omega_3} \int_{\mathbb{R}^4} \frac{w^{p_2}_{\lam}(y) e^{q_2 z_{\lam}(y)}}{|x-y|^2} dy + \frac{2(\mu - \alpha)}{|x|^2} \\
& = \frac{1}{2 \omega_3} \int_{B_{\lam}} \frac{w^{p_2}(y) e^{q_2 z(y)}}{|x-y|^2} dy + \frac{1}{2 \omega_3} \int_{B_{\lam}^c} \frac{w^{p_2}(y) e^{q_2 z(y)}}{|x-y|^2} dy \\
& - \frac{1}{2 \omega_3} \int_{B_{\lam}} \frac{w^{p_2}_{\lam}(y) e^{q_2 z_{\lam}(y)}}{|x-y|^2} dy - \frac{1}{2 \omega_3} \int_{B_{\lam}^c} \frac{w^{p_2}_{\lam}(y) e^{q_2 z_{\lam}(y)}}{|x-y|^2} dy + \frac{2(\mu - \alpha)}{|x|^2} \\
& = \frac{1}{2 \omega_3} \int_{B_{\lam}} \left[\frac{1}{|x-y|^2} - \frac{1}{\left|x-\frac{\lam^2 y}{|y|^2}\right|^2}\right][w^{p_2}(y) e^{q_2 z(y)} - w^{p_2}_{\lam}(y) e^{q_2 z_{\lam}(y)}] dy + \frac{2(\mu - \alpha)}{|x|^2}
\end{aligned}
\end{equation}
and
\[
\Delta^2 S_\lam^z(x) = w^{p_2}(x) e^{q_2 z(x)} - w^{p_2}_{\lam}(x) e^{q_2 z_{\lam}(x)}.
\]

\begin{pro}\label{b40}
Assume that $ \alpha < \mu$. Then for $ \lam > 0$ large enough, we have
\[
S_\lam^w(x) \ge 0, \ S_\lam^z(x) \ge 0 \ \text{and} \ -\Delta S_\lam^z(x) \ge 0 \ \text{in} \ B_\lam(0) \setminus \{0\}.
\]
\end{pro}
\noindent \emph{Proof.}

Recall the following lemma from \cite{CLn,Yu}(see Lemma 2.1 of \cite{CLn} and Lemma 3.3 of \cite{Yu}).

\begin{lemma}\label{b38}
Suppose that $ w(x)$ satisfies
\begin{equation}
\left\{
\begin{aligned}
& -\Delta w(x) \geq 0, \quad \text{in} \  \mathbb{R}^4 \setminus \{0\}\\
& w(x) > 0 , \quad \text{in} \ \mathbb{R}^4 \setminus \{0\}.
\nonumber
\end{aligned}
\right.
\end{equation}
If we denote $ \eps = \min\limits_{\ptl B_R}w(x)$, then $ w(x)\geq \eps$ for $ x\in B_R \setminus \{0\}$.
\end{lemma}

Step 1: There exists a $ R_0 > 0$, such that $ S_\lam^w(x) \ge 0$, $ S_\lam^z(x) \ge 0$ and $ -\Delta S_\lam^z(x) \ge0$ for $ R_0 \leq |x| \leq \frac{\lam}{2}$  and $\lambda$ large enough.

In fact,
\begin{align*}
S_\lam^w(x) & = \frac{1}{|x|^2} u \left(\frac{x}{|x|^2}\right) - \frac{1}{\lam^2} u\left(\frac{x}{\lam^2}\right) \\
& = \left(\frac{1}{|x|^2} -\frac{1}{\lam^2} \right) u \left(\frac{x}{|x|^2}\right) + \frac{1}{\lam^2}\left[u \left(\frac{x}{|x|^2}\right) - u\left(\frac{x}{\lam^2}\right) \right] \\
& \geq \frac{3}{\lam^2} [u(0)+ o(1)] - \frac{C}{\lam^2}\left(\frac{1}{R_0} - \frac{R_0}{\lam^2} \right) \\
& \geq \frac{3}{\lam^2} [u(0)+ o(1)] - \frac{C}{\lam^2 R_0} > 0
\end{align*}
for $ \lam$ and $ R_0$ large enough, where $ C= \max_{B_1} |\nabla u|$. Similarly, we find
\begin{align*}
S_\lam^z(x) & =v \left(\frac{x}{|x|^2}\right) - \mu \ln|x| - v \left( \frac{x}{\lam^2}\right) + \mu \ln \lam \\
& = v \left(\frac{x}{|x|^2}\right)- v \left( \frac{x}{\lam^2}\right) + \mu \ln \left( \frac{\lam}{|x|} \right)\\
& \geq \mu \ln 2 - \frac{C}{R_0} \\
& > 0.
\end{align*}
Finally, one gets
\begin{align*}
- \Delta S_\lam^z(x) & =- \Delta z(x) + \Delta z_\lam(x) \\
& = \frac{4x}{|x|^4}\nabla v \left(\frac{x}{|x|^2}\right) - \frac{1}{|x|^4}\Delta v \left(\frac{x}{|x|^2} \right) + \frac{2\mu}{|x|^2} + \frac{1}{\lam^4} \Delta v\left(\frac{x}{\lam^2} \right)\\
& \geq \frac{2\mu}{|x|^2} -\frac{4C}{|x|^3} -\frac{C}{|x|^4} -\frac{C}{\lam^4} \\
& > 0
\end{align*}
for $ \lam$ and $ R_0$ large enough, where $ C= \max_{B_1} \{ |\nabla v|, |\Delta v|\}$. This completes the proof of Step 1.

Step 2: $ S_\lam^w(x) \ge 0$, $ S_\lam^z(x) \ge 0$ and $ -\Delta S_\lam^z(x) \ge 0$ for $ \frac{\lam}{2} \le |x| \le \lam$.

In fact, we can get $ -\Delta S_\lam^z(x) \ge 0$ in the same way as Step 1. Since $ S_\lam^z$ satisfies
\begin{equation}
\left\{
\begin{aligned}
&-\Delta S_\lam^z \ge 0 \quad \text{in} \ B_\lam \setminus B_{\frac{\lam}{2}} \\
&S_\lam^z = 0 \qquad \text{on} \ \ptl B_\lam, \\
&S_\lam^z \ge 0 \qquad \text{on} \ptl B_{\frac{\lam}{2}},
\nonumber
\end{aligned}
\right.
\end{equation}
then we infer from the maximum principle that $ S_\lam^z \ge 0$ for $  \frac{\lam}{2} \leq |x| \leq \lam$.

Finally, since $ S_\lam^w(x)$ satisfies
\begin{align*}
-\Delta S_\lam^w(x) &  =  w^{p_1} e^{q_1 z} - w_{\lam}^{p_1} e^{q_1 z_{\lam}} \\
& = w^{p_1}( e^{q_1 z} - e^{q_1 z_{\lam}}) + (w^{p_1} - w_{\lam}^{p_1}) e^{q_1 z_{\lam}} \\
& = q_1 w^{p_1} e^{q_1 \xi_z}(z- z_{\lam}) + p_1 \xi_w^{p_1 - 1} (w-w_{\lam}) e^{q_1 z_{\lam}},
\end{align*}
where $ \xi_z$ is between $ z(x)$ and $ z_{\lam}(x)$ and $ \xi_w$ is between $ w(x)$ and $ w_{\lam}(x)$. We see that $ q_1 w^{p_1} e^{q_1 \xi_z}(z- z_{\lam}) \ge 0$ and $ p_1 \xi_w^{p_1 - 1}e^{q_1 z_{\lam}} \geq 0$ in $ B_\lam \setminus B_{\frac{\lam}{2}}$. We obatin
\begin{equation}
\left\{
\begin{aligned}
& - \Delta S_\lam^w(x) - c(x) S_\lam^w(x) \ge 0 \quad \text{in}  \ B_\lam \setminus B_{\frac{\lam}{2}} \\
&  S_\lam^w(x) = 0 \qquad \text{on} \ \ptl B_\lam \\
&S_\lam^w \ge 0 \qquad \text{on} \ \ptl B_{\frac{\lam}{2}},
\nonumber
\end{aligned}
\right.
\end{equation}
where $ c(x) = p_1 \xi_w^{p_1 - 1}e^{q_1 z_{\lam}} \geq 0$ is bounded. By a simple calculation, we observe that $ c(x) = O(\frac{1}{\lam^4})$, $ x \in  B_\lam \setminus B_{\frac{\lam}{2}} $.

Let $ \Omega = B_{\frac{5}{4}} \setminus B_{\frac{1}{4}}$. Let $ \lam_1$ be the first Dirichlet eigenvalue of $ - \Delta$ in $ \Omega$ and $ \varphi(x)$ be an eigenfunction corresponding to $ \lam_1$.
We also assume $\varphi>0$ in $\Omega$ and $\|\varphi\|_{L^\infty(\Omega)}=1$.
Let $ \psi(x) = \varphi (\frac{x}{\lam})$ and $ \Omega_{\lam} = B_{\frac{5}{4}\lam} \setminus B_{\frac{1}{4}\lam}$. Then, $ \psi(x)$ satisfies
\begin{equation}
\left\{
\begin{aligned}
&- \Delta \psi (x)=\frac{\lam_1}{\lam^2} \psi(x) ,  \quad \text{in} \ \Omega_{\lam} \\
& \psi (x)  = 0, \qquad \text{on} \ \ptl \Omega_{\lam}.
\nonumber
\end{aligned}
\right.
\end{equation}
Let $S(x)=\frac{S_\lam^w(x)}{\psi(x)}$. Then $ S(x)$ solves
\begin{equation*}
\begin{cases}
&-\Delta S-2\nabla S\cdot \nabla \ln\psi+\left(\frac{\lambda_1}{\lambda^2}-c(x)\right)S \ge 0,\quad \text{in}\quad B_\lambda\backslash B_{\frac{\lambda}{2}},\\
& S \ge 0,\quad \text{on}\quad \ptl B_{\frac{\lambda}{2}}, \\
& S = 0,\quad \text{on}\quad \ptl B_\lam.
\end{cases}
\end{equation*}
Notice that $\frac{\lambda_1}{\lambda^2}-c(x)\ge 0$ for $\lambda$ large enough. By standard maximum principle, one gets
$S(x)\ge 0$, i.e. $S_\lambda^w(x)\ge 0$ in $B_\lambda\backslash B_{\frac{\lambda}{2}}$.

Step 3: $ S_\lam^w(x) \ge 0$, $ S_\lam^z(x) \ge 0$ and $ -\Delta S_\lam^z(x) \ge 0$ in $ B_{R_0} \setminus \{0\}$.

Claim:  there exists a $ R_2$ small enough, such that $ S_\lam^z(x) \ge 0$ for $ 0< |x| \leq R_2$.
\\ Since $$v\left(\frac{x}{|x|^2}\right)=[\alpha + o(1)] \ln |x|,\quad 0< |x| \leq R_2,$$
then
\begin{equation*}
\begin{split}
S_\lam^z(x)=&v\left(\frac{x}{|x|^2}\right) -v\left(\frac{x}{\lam^2}\right) - \mu \ln |x| + \mu \ln \lambda
\\ =&[\alpha - \mu + o(1)]\ln |x| -v\left(\frac{x}{\lam^2}\right)  + \mu \ln \lam >0
\end{split}
\end{equation*}
for $ 0 < |x| \leq R_2$ small enough. This proves the claim.

Fixing $R_2>0$ as above and choosing $\lambda>0$ large enough, one gets from the definition of $S_\lambda^z$ such that $S_\lambda^z(x) \ge 0$, $\forall R_2\le |x|\le R_0$.

Since
\[
- \Delta \left(\frac{1}{|x|^2}u\left(\frac{x}{|x|^2}\right)\right) = -\frac{1}{|x|^6} \Delta u\left(\frac{x}{|x|^2}\right) = \frac{u^{p_1}\left(\frac{x}{|x|^2}\right) e^{q_1 v\left(\frac{x}{|x|^2}\right)}}{|x|^6} > 0
\]
for $ x\neq 0$, it follows from Lemma \ref{b38} that
\[
\frac{1}{|x|^2}u\left(\frac{x}{|x|^2}\right) \geq \eps > 0 \quad \text{in} \ B_{R_0}\setminus \{0\},
\]
where $ \eps = \min\limits_{\ptl B_{R_0}}\frac{1}{|x|^2}u\left(\frac{x}{|x|^2}\right) >0$. So we conclude that
\[
S_\lam^w(x) = \frac{1}{|x|^2}u\left(\frac{x}{|x|^2}\right) - \frac{1}{\lam^2}u\left(\frac{x}{\lam^2}\right) \geq \eps - \frac{1}{\lam^2}u\left(\frac{x}{\lam^2}\right) > 0
\]
in $ B_{R_0}\setminus \{0\}$ for $ \lam$ large enough.

Similarly, since
\[
|x|^2 (- \Delta S_\lam^z(x) ) = |x|^2 \left[ \frac{4x}{|x|^4}\nabla v \left(\frac{x}{|x|^2}\right) - \frac{1}{|x|^4}\Delta v \left(\frac{x}{|x|^2} \right) + \frac{2\mu}{|x|^2} + \frac{1}{\lam^4} \Delta v\left(\frac{x}{\lam^2} \right) \right],
\]
then we infer from Lemma \ref{b39} that
\[
|x|^2 (- \Delta S_\lam^z(x) )  \geq  2(\mu - \alpha)  + \frac{|x|^2}{\lam^4} \Delta v\left(\frac{x}{\lam^2} \right) > 0
\]
for $ \lam$ large enough. Hence we deduce that
\[
- \Delta S_\lam^z(x) \ge 0 \quad \text{in} \quad B_{R_0}\setminus \{0\}.
\]
Proposition \ref{b40} follows from Step 1 to Step 3. \hfill $ \Box$

Next, we give a technical lemma.

\begin{lemma}\label{b41} (see Lemma 11.2 of \cite{LiyZhang})
(1)Suppose $ u \in C^1(\mathbb{R}^4)$, if for all $ b\in \mathbb{R}^4$ and $ \lam > 0$, the following inequality holds
\[
\frac{1}{|x|^2} u_b \left(\frac{x}{|x|^2}\right) - \frac{1}{\lam^2} u_b \left(\frac{x}{\lam^2}\right) \geq 0, \ \forall x \in B_\lam \setminus \{0\},
\]
then we have $ u(x) \equiv C$, where $ u_b(x)=u(x+b)$.

(2) Suppose $ v \in C^1(\mathbb{R}^4)$, if for all $ b\in \mathbb{R}^4$ and $ \lam > 0$, the following inequality holds
\[
v_b \left(\frac{x}{|x|^2}\right)-\mu \ln |x| - v_b \left(\frac{x}{\lam^2}\right) + \mu \ln \lam \geq 0, \quad x \in B_\lam \setminus \{0\}
\]
then we have $ v(x) \equiv C$, where $ v_b(x)=v(x+b)$.
\end{lemma}

Now we define
\begin{equation*}
\begin{split}
&  u_b(x) = u(x+b),\quad v_b(x) = v(x+b)\\
&  w_b(x) = \frac{1}{|x|^2} u_b \left(\frac{x}{|x|^2}\right),\quad z_b(x) = v_b \left(\frac{x}{|x|^2}\right) - \mu \ln|x|,\\
& w_{\lam,b}(x) = \frac{\lam^2}{|x|^2} w_b\left(\frac{\lam^2 x}{|x|^2}\right) = \frac{1}{\lam^2} u_b\left(\frac{x}{\lam^2}\right),\\
& z_{\lam,b}(x) = z_b\left(\frac{\lam^2 x}{|x|^2}\right) - \mu \ln \frac{|x|}{\lam} = v_b\left( \frac{x}{\lam^2}\right) - \mu \ln \lam,\\
& S_{\lam,b}^w(x) = w_b(x) - w_{\lam,b}(x),\quad S_{\lam,b}^z(x) = z_b(x) - z_{\lam,b}(x).
\end{split}
\end{equation*}

For fixed $ b\in \mathbb{R}^4$, we define
\[
\lam_b = \inf\{\lam > 0 \ | \ S_{\mu,b}^w(x)\ge 0,  S_{\mu,b}^z(x)\ge 0, -\Delta  S_{\mu,b}^z(x)\ge0, \ \text{in} \ B_\mu \setminus \{0\}, \lam \leq \mu < \infty\}.
\]

\begin{pro}\label{b55}
There exists a vector $ b \in \mathbb{R}^4$, such that $ \lam_b > 0$.
\end{pro}

\noindent \emph{Proof.} We prove it by contradiction. Suppose on the contrary, then for any $ b\in \mathbb{R}^4$, we have $ \lam_b =0$. By the definition of $ \lam_b$, we get
\[
S_{\lam,b}^w(x)\ge 0 \ \text{and} \ S_{\lam,b}^z(x)\ge 0
\]
for any $ \lam > 0$ and $ x\in B_\lam \setminus \{0\}$. Then we infer from Lemma \ref{b41} that $ u(x)\equiv C_1$ and $ v(x)\equiv C_2$. This contradicts $ \int_{\mathbb{R}^4}  u^{p_2}(x) e^{q_2 v(x)} dx < \infty$ unless $ C_1 = 0$. \hfill $ \Box$

\begin{pro}\label{b46}
If $ (u,v)$ is a nontrivial solution of problem \eqref{b1}, then $ \alpha < \mu$ can't occur.
\end{pro}

\noindent \emph{Proof.} Suppose on the contrary that $ \alpha < \mu$, then it follows from Proposition \ref{b55} that there exists a vector $ b \in \mathbb{R}^4$, such that $ \lam_b > 0$. Without loss of generality, we assume $b=0$. Then
\be\label{b43}
S_{\lam_0}^w(x) \geq 0, \ \ S_{\lam_0}^z(x)\geq 0 \ \text{and} \  -\Delta  S_{\lam_0}^z(x)\geq 0
\ee
in $ B_{\lam_0} \setminus \{0\} $. Moreover, we deduce from equation \eqref{b42} and \eqref{b43} that
\[
- \Delta S_{\lam_0}^z(x) > 0 \quad \text{in}\quad B_{\lam_0} \setminus \{0\}.
\]
Then the maximum principles implies
\[
S_{\lam_0}^z(x) > 0,\quad \text{in}\quad B_{\lam_0} \setminus \{0\}.
\]
By \eqref{b43}, we conclude that
\[
-\Delta S_{\lam_0}^w(x) \geq 0,\quad \text{in}\quad B_{\lam_0} \setminus \{0\}.
\]
This implies
\[
S_{\lam_0}^w(x) > 0 \quad \text{or} \quad  S_{\lam_0}^w(x) \equiv 0 \quad \text{in}\quad  B_{\lam_0} \setminus \{0\}.
\]
If $ S_{\lam_0}^w(x) \equiv 0$, it is clear that
\[
-\Delta S_{\lam_0}^w(x) = w^{p_1}(x) (e^{q_1z(x)} - e^{q_1 z_{\lam_0}(x)})> 0,\quad \text{in}\quad B_{\lam_0} \setminus \{0\}.
\]
This leads to a contradiction. Therefore, we obtain $ S_{\lam_0}^w(x) > 0$ in $ B_{\lam_0} \setminus \{0\} $.
By the definition of $ \lam_0$, one of the following three cases may occur.

(i) $ \exists \lam_k < \lam_0$, $ \lam_k \rightarrow \lam_0$ with $ \inf\limits_{B_{\lam_k}\setminus \{0\}} S_{\lam_k}^w(x) < 0$.

(ii) $ \exists \lam_k < \lam_0$, $ \lam_k \rightarrow \lam_0$ with $ \inf\limits_{B_{\lam_k}\setminus \{0\}} S_{\lam_k}^z(x) < 0$.

(iii) $ \exists \lam_k < \lam_0$, $ \lam_k \rightarrow \lam_0$, such that $ S_{\lam_k}^w(x) \geq 0$, $ S_{\lam_k}^z(x) \geq 0$ for $ x \in B_{\lam_k}\setminus \{0\}$, but $ \inf\limits_{B_{\lam_k}\setminus \{0\}} - \Delta S_{\lam_k}^z(x) < 0$. \\
Now we show that each case will lead to a contradiction.

If case (i) occur, then we deduce from $ S_{\lam_0}^w(x) > 0 $ in $ B_{\lam_0} \setminus \{0\} $ and the Hopf Lemma that
\be\label{b44}
\frac{\ptl S_{\lam_0}^w}{\ptl \nu} (x) < 0
\ee
on $ \ptl B_{\lam_0}$, where $ \nu$ is the unit outer normal direction.
\par Claim:  there exists a $ \gamma= \gamma(\frac{\lam_0}{2}) >0$, such that $ S_{\lam_k}^w(x) \geq \frac{\gamma}{2}$, $ \forall x\in B_{\frac{\lam_0}{2}} \setminus \{0\}$. \\
Let
\[
\gamma = \min_{\ptl B_{\frac{\lam_0}{2}}} S_{\lam_0}^w(x) > 0.
\]
We define
\[
h(x) = \gamma - \frac{r^2}{|x|^2}\gamma \quad \text{in} \quad B_{\frac{\lam_0}{2}} \setminus B_r
\]
with $ r$ small. Then, $ k(x) = S_{\lam_0}^w(x) - h(x)$ satisfies
\begin{equation}
\left\{
\begin{aligned}
& - \Delta k(x) = - \Delta S_{\lam_0}^w(x) \geq 0 \quad \text{in} \quad B_{\frac{\lam_0}{2}} \setminus B_r, \\
& k(x) = S_{\lam_0}^w(x) > 0 \qquad \text{on} \quad \ptl B_r \\
& k(x)> 0 \qquad \qquad \text{on} \quad \ptl B_{\frac{\lam_0}{2}}.
\nonumber
\end{aligned}
\right.
\end{equation}
Hence, by the maximum principle and letting $ r\rightarrow 0$, one gets $S_{\lam_0}^w(x) \ge \gamma$, in $B_{\frac{\lambda_0}{2}}\backslash\{0\}$.  Then
\begin{equation*}
\begin{split}
S_{\lambda_k}^w(x)=&w(x)-w_{\lambda_k}(x)\\
=&S^{w}_{\lambda_0}(x)+w_{\lambda_0}(x)-w_{\lambda_k}(x)\\
\ge & \frac{\gamma}{2}
\end{split}
\end{equation*}
provided $\lambda_k$ is close enough to $\lambda_0$. This proves the claim.

On the other hand, since $ \inf\limits_{B_{\lam_k}\setminus \{0\}} S_{\lam_k}^w(x) < 0$, then we infer from the claim that there exists an $ x_k \in B_{\lam_k}\setminus B_{\frac{\lam_0}{2}} $ such that
\[
S_{\lam_k}^w(x_k) = \inf_{B_{\lam_k}\setminus \{0\}} S_{\lam_k}^w(x) < 0.
\]
In particular, we have $ \nabla S_{\lam_k}^w(x_k) = 0$. We assume that, up to a subsequence, $ x_k \rightarrow \bar{x}$, then we obtain $ \nabla S_{\lam_0}^w(\bar{x}) =0 $ and $ S_{\lam_0}^w(\bar{x}) =0$. Hence $ \bar{x} \in \ptl B_{\lam_0}$. However, this contradicts equation \eqref{b44}. Hence, case (i) can't occur.

Next, we show that case (ii) can not occur. We first claim that there exists a $ r_0 >0$ such that
\[
S_{\lam_0}^z(x) > 1 \quad \text{in} \quad B_{r_0} \setminus \{0\}.
\]
Since
\[
\lim_{|x|\rightarrow 0}\frac{v\left(\frac{x}{|x|^2}\right)}{\ln |x|}=\alpha,
\]
 then for small $ r_0>0$, we get
\begin{align*}
v\left(\frac{x}{|x|^2}\right) - \mu \ln|x| & = (\alpha + o(1) - \mu) \ln |x| \\
& \geq (\mu - \alpha + o(1)) (-\ln r_0) \\
& \geq \frac{\alpha - \mu}{2} \ln r_0
\end{align*}
for $ 0 < |x| < r_0$. Set $ C= \sup\limits_{B_{\frac{1}{\lam_0}}} v$. Then we can further choose $ r_0$ small enough, such that
\begin{equation}\label{b45}
\begin{aligned}
S_{\lam_0}^z(x) & = v \left(\frac{x}{|x|^2}\right) - \mu \ln|x| - v\left( \frac{x}{\lam_0^2}\right) + \mu \ln \lam_0 \\
& \geq \frac{\alpha - \mu}{2} \ln r_0 -C + \mu \ln \lam_0 \\
& > 1
\end{aligned}
\end{equation}
for $ 0<|x| < r_0$. This proves the claim.

If case (ii) occur, then there exist some $ \lam_k < \lam_0$, $ \lam_k \rightarrow \lam_0$ such that
\[
\inf\limits_{B_{\lam_k}\setminus \{0\}} S_{\lam_k}^z(x) < 0.
\]
By \eqref{b45} and similar argument as case (i), we have $ S_{\lam_k}^z(x) \geq \frac{1}{2}$ for $ 0 < |x| \leq r_0$ and $ k$ large enough. Hence $ \inf\limits_{B_{\lam_k}\setminus \{0\}} S_{\lam_k}^z(x)$ is attained at some $ x_k \in B_{\lam_k} \setminus B_{r_0}$. Therefore,
\[
\nabla S_{\lam_k}^z(x_k) = 0.
\]
We assume that, up to a subsequence, $ x_k \rightarrow \bar{x}$ as $ k \rightarrow \infty$, then
\[
\nabla S_{\lam_0}^z(\bar{x}) = 0 \quad \text{and} \quad  S_{\lam_0}^z(\bar{x}) = 0,
\]
which implies that $ \bar{x} \in \ptl B_{\lam_0}$. But this contradicts the Hopf Lemma.

Finally, we show that case (iii) can't occur. For $ 0< |x| < \lam_k$, it is clear that
\begin{align*}
- \Delta S_{\lam_k}^z &= \frac{1}{2 \omega_3} \int_{B_{\lam_k}} \left[\frac{1}{|x-y|^2} - \frac{1}{\left|x-\frac{\lam_k^2 y}{|y|^2}\right|^2}\right][w^{p_2}(y) e^{q_2 z(y)} - w^{p_2}_{\lam_k}(y) e^{q_2 z_{\lam_k}(y)}] dy + \frac{2(\mu - \alpha)}{|x|^2} \\
& > 0,
\end{align*}
which contradicts
\[
\inf\limits_{B_{\lam_k}\setminus \{0\}} - \Delta S_{\lam_k}^z(x) < 0.
\]
\hfill $ \Box$

\section{Fast decay is impossible}
This section is devoted to exclude the case $ \alpha > \mu$. We prove the conclusion by contradiction. Therefore, we assume that $ \alpha > \mu$ in this section.

\begin{lemma}\label{b48}
Let
\[
\beta = \frac{1}{2 \omega_3} \int_{\mathbb{R}^4}  u^{p_1}(x) e^{q_1 v(x)} dx.
\]
Assume that $ \alpha > \mu$, then we have
\[
\lim\limits_{|x| \rightarrow \infty} (|x|^2 u(x) - \beta ) = 0.
\]
\end{lemma}
\noindent \emph{Proof.} By Lemma \ref{b14}, we can write
\[
|x|^2 u(x) - \beta = \frac{1}{2 \omega_3} \int_{\mathbb{R}^4} \frac{|x|^2 - |x-y|^2}{|x-y|^2} u^{p_1}(y) e^{q_1 v(y)} dy.
\]
Consider $ M$ large enough, and set
\[
D_1 = \left \{y \in \mathbb{R}^4 | |x-y| \leq \frac{|x|}{2} \right\}, \quad D_2 = \left \{y \in \mathbb{R}^4 | |x-y| > \frac{|x|}{2},|y|\leq M \right\}
\]
and
\[
D_3 = \left\{y \in \mathbb{R}^4 | |x-y| > \frac{|x|}{2},|y|> M \right\}.
\]

For $ y \in D_1$, we have $ \frac{|x|}{2} \leq |y| \leq \frac{3}{2}|x|$.
We can choose $ \eps$ small enough such that $ \alpha - \eps > \mu $.
By Lemma \ref{b22} and Lemma \ref{b15}, we find
\begin{align*}
& \int_{D_1} \frac{|x|^2 - |x-y|^2}{|x-y|^2} u^{p_1}(y) e^{q_1 v(y)} dy \\
& \leq  \int_{D_1} u^{p_1}(y) e^{q_1 v(y)} dy + |x|^2  \int_{D_1} \frac{u^{p_1}(y) e^{q_1 v(y)} }{|x-y|^2} dy \\
& \leq o(1) + C|x|^{2  -\mu q_1} \int_{D_1} \frac{1}{|x-y|^2}dy \\
& \leq o(1) + C|x|^{4  -\mu q_1} \rightarrow 0
\end{align*}
as $ |x| \rightarrow \infty$. In getting the above inequality, we also used $ \mu = \frac{6-2p_1}{q_1}$, $ p_1 < 1$ and $ v(x) \sim -\alpha \ln |x|$ at $ \infty $.

For $ y \in D_2$, one has
\[
 \frac{|x|^2 - |x-y|^2}{|x-y|^2} \rightarrow 0
\]
and
\[
\int_{D_2} \frac{|x|^2 - |x-y|^2}{|x-y|^2} u^{p_1}(y) e^{q_1 v(y)} dy \rightarrow 0
\]
as $ |x| \rightarrow \infty$.

For $ y \in D_3$, we have $ \frac{|x|}{|x-y|} < 2$, which further implies
\[
\left|\frac{|x|^2 - |x-y|^2}{|x-y|^2}\right| < 3.
\]
Hence, we choose $ M$ large enough and we deduce that
\[
\int_{D_3} \frac{|x|^2 - |x-y|^2}{|x-y|^2} u^{p_1}(y) e^{q_1 v(y)} dy = o(1).
\]
This completes the proof of this lemma. \hfill $ \Box$

\begin{lemma}\label{b51}
(i) $ x\cdot \nabla v(x) + \alpha \rightarrow 0 $ as $ |x| \rightarrow \infty$. \\
(ii) $ |x|^2 \Delta v(x) + 2 \alpha \rightarrow 0 $ as $ |x| \rightarrow \infty$.
\end{lemma}

\noindent \emph{Proof.} By \eqref{b28}, we can write
\[
x\cdot \nabla v(x) + \alpha = \frac{1}{4 \omega_3} \int_{\mathbb{R}^4} \frac{|x-y|^2 - |x|^2 + x\cdot y}{|x-y|^2} u^{p_2}(y) e^{q_2 v(y)} dy.
\]
We divide the integral domain $ \mathbb{R}^4$ into $ D_1 \cup D_2 \cup D_3$ as Lemma \ref{b48}.

For $ y \in D_1$, we obtain
\begin{align*}
& \int_{D_1} \frac{|x-y|^2 - \langle x,x-y \rangle}{|x-y|^2} u^{p_2}(y) e^{q_2 v(y)} dy \\
& \leq  \int_{D_1} u^{p_2}(y) e^{q_2 v(y)} dy + |x|  \int_{D_1} \frac{u^{p_2}(y) e^{q_2 v(y)} }{|x-y|} dy \\
& \leq o(1) + C|x|^{-7} \int_{D_1} \frac{1}{|x-y|}dy \\
& \leq o(1) + C|x|^{- 4} \rightarrow 0,
\end{align*}
as $ |x| \rightarrow \infty$. In order to get the above inequality, we also used $ u(x) \sim \frac{\beta}{|x|^2}$ at $ \infty$. For $ y \in D_2 \cup D_3$, the proof is similar to Lemma \ref{b48}.

By \eqref{b18}, it is clear that
\[
|x|^2 \Delta v(x) + 2 \alpha = \frac{1}{2 \omega_3} \int_{\mathbb{R}^4} \frac{|x-y|^2 - |x|^2 }{|x-y|^2} u^{p_2}(y) e^{q_2 v(y)} dy.
\]
In a similar way, we can prove
\[
\int_{D_i} \frac{|x-y|^2 - |x|^2 }{|x-y|^2} u^{p_2}(y) e^{q_2 v(y)} dy \rightarrow 0, \quad i=1,2,3
\]
as $ |x| \rightarrow \infty$. \hfill $ \Box$

Define
\[
\widetilde{S}^u_{\lam} := u(x) - u_{\lam}(x) =u(x) -  \frac{\lam^2}{|x|^2} u \left(\frac{\lam^2 x}{|x|^2} \right)
\]
\[
\widetilde{S}^v_{\lam} : = v(x) - v_{\lam}(x) = v(x) - v\left(\frac{\lam^2 x}{|x|^2} \right) + \mu \ln|x| - \mu \ln \lam,
\]
then we have the following result.

\begin{lemma}\label{b52} We have \\
(i) $ \Delta \left( v \left(\frac{\lam^2 x}{|x|^2} \right)\right)= \frac{2 \alpha}{|x|^2} - \frac{1}{2 \omega_3 } \int_{\mathbb{R}^4} \frac{u^{p_2}(y) e^{q_2 v(y)}}{\left|x-\frac{\lam^2 y }{|y|^2}\right|^2} dy$  \\
(ii) $ - \Delta \widetilde{S}^v_{\lam} = \frac{1}{2 \omega_3 } \int_{B_{\lam}} \left[\frac{1}{|x-y|^2} - \frac{1}{\left|x-\frac{\lam^2 y }{|y|^2}\right|^2}\right] [u^{p_2}(y) e^{q_2 v(y)} - u_{\lam}^{p_2} e^{q_2 v_{\lam}} ] dy + \frac{2(\alpha-\mu)}{|x|^2}$ \\
for $ x\in \mathbb{R}^4 \setminus \{0\}$.
\end{lemma}

\noindent \emph{Proof.}
A direct calculation shows that
\begin{align*}
\Delta \left( v \left(\frac{\lam^2 x}{|x|^2} \right) \right) & = \frac{\lam^4}{|x|^4} \Delta v \left(\frac{\lam^2 x}{|x|^2} \right)  - 4\lam^2 \frac{x \cdot \nabla v \left(\frac{\lam^2 x}{|x|^2} \right)}{|x|^4} \\
& = - \frac{\lam^4}{2 \omega_3 |x|^4} \int_{\mathbb{R}^4} \frac{u^{p_2}(y) e^{q_2 v(y)}}{\left|\frac{\lam^2 x}{|x|^2} -y \right|^2} dy + \sum^4_{i=1} \frac{\lam^2 x_i}{\omega_3 |x|^4 } \int_{\mathbb{R}^4} \frac{\frac{\lam^2 x_i}{|x|^2}-y_i}{\left|\frac{\lam^2 x}{|x|^2}-y\right|^2} u^{p_2}(y) e^{q_2 v(y)} dy \\
& = \frac{\lam^4}{2 \omega_3 |x|^4} \int_{\mathbb{R}^4} \frac{u^{p_2}(y) e^{q_2 v(y)}}{\left|\frac{\lam^2 x}{|x|^2} -y\right|^2} dy - \frac{\lam^2 }{\omega_3 |x|^4 } \int_{\mathbb{R}^4} \frac{x \cdot y}{\left|\frac{\lam^2 x}{|x|^2}-y\right|^2} u^{p_2}(y) e^{q_2 v(y)} dy \\
& = \frac{1}{2 \omega_3 } \int_{\mathbb{R}^4} \frac{u^{p_2}(y) e^{q_2 v(y)}}{\left|x - \frac{|x|^2}{\lam^2} y\right|^2} dy - \frac{1}{\omega_3 } \int_{\mathbb{R}^4} \frac{ \frac{1}{\lam^2} x \cdot y}{\left|x - \frac{|x|^2}{\lam^2 } y \right|^2} u^{p_2}(y) e^{q_2 v(y)} dy \\
& =  \frac{1}{2 \omega_3 } \int_{\mathbb{R}^4} \frac{1- \frac{2x\cdot y}{\lam^2}}{\left|x - \frac{|x|^2}{\lam^2 } y \right|^2} u^{p_2}(y) e^{q_2 v(y)} dy \\
& =  \frac{1}{2 \omega_3 } \int_{\mathbb{R}^4} \frac{\frac{1}{\lam^4} \left[|x|^2 \left|\frac{\lam^2 x}{|x|^2} - y \right|^2 - |x|^2 |y|^2\right]}{\left|x - \frac{|x|^2}{\lam^2 } y \right|^2} u^{p_2}(y) e^{q_2 v(y)} dy \\
& =  \frac{1}{2 \omega_3 } \int_{\mathbb{R}^4} \left[\frac{1}{|x|^2} - \frac{|x|^2 |y|^2}{\lam^4 \left|x - \frac{|x|^2}{\lam^2 } y \right|^2}\right] u^{p_2}(y) e^{q_2 v(y)} dy \\
& = \frac{2 \alpha}{|x|^2} - \frac{1}{2 \omega_3 } \int_{\mathbb{R}^4} \frac{u^{p_2}(y) e^{q_2 v(y)}}{\left|x-\frac{\lam^2 y }{|y|^2} \right|^2} dy
\end{align*}
and
\begin{align*}
- \Delta \widetilde{S}^v_{\lam} & = - \Delta v(x) + \Delta \left( v \left(\frac{\lam^2 x}{|x|^2} \right)\right) - \frac{2\mu}{|x|^2}  \\
& = \frac{1}{2 \omega_3} \int_{\mathbb{R}^4} \frac{u^{p_2}(y) e^{q_2 v(y)}}{|x-y|^2} dy - \frac{1}{2 \omega_3 } \int_{\mathbb{R}^4} \frac{u^{p_2}(y) e^{q_2 v(y)}}{\left|x-\frac{\lam^2 y }{|y|^2} \right|^2} dy + \frac{2(\alpha-\mu)}{|x|^2} \\
& = \frac{1}{2 \omega_3} \int_{B_{\lam}} \frac{u^{p_2}(y) e^{q_2 v(y)}}{|x-y|^2} dy +\frac{1}{2 \omega_3} \int_{B_{\lam}^c} \frac{u^{p_2}(y) e^{q_2 v(y)}}{|x-y|^2} dy\\
& - \frac{1}{2 \omega_3 } \int_{B_{\lam}} \frac{u^{p_2}(y) e^{q_2 v(y)}}{\left|x-\frac{\lam^2 y }{|y|^2}\right|^2} dy - \frac{1}{2 \omega_3 } \int_{B^c_{\lam}} \frac{u^{p_2}(y) e^{q_2 v(y)}}{\left|x-\frac{\lam^2 y }{|y|^2} \right|^2} dy + \frac{2(\alpha-\mu)}{|x|^2} \\
& = \frac{1}{2 \omega_3} \int_{B_{\lam}} \frac{u^{p_2}(y) e^{q_2 v(y)}}{|x-y|^2} dy +\frac{1}{2 \omega_3} \int_{B_{\lam}} \frac{u^{p_2}\left(\frac{\lam^2 z}{|z|^2}\right) e^{q_2 v\left(\frac{\lam^2 z}{|z|^2}\right)}}{\left|x-\frac{\lam^2 z}{|z|^2}\right|^2} \frac{\lam^8}{|z|^8} dz \\
& - \frac{1}{2 \omega_3 } \int_{B_{\lam}} \frac{u^{p_2}(y) e^{q_2 v(y)}}{\left|x-\frac{\lam^2 y }{|y|^2}\right|^2} dy - \frac{1}{2 \omega_3 } \int_{B_{\lam}}\frac{u^{p_2}\left(\frac{\lam^2 z}{|z|^2} \right) e^{q_2 v\left(\frac{\lam^2 z}{|z|^2} \right)}}{|x-z|^2}\frac{\lam^8}{|z|^8} dz + \frac{2(\alpha-\mu)}{|x|^2} \\
& = \frac{1}{2 \omega_3 } \int_{B_{\lam}}\left [\frac{1}{|x-y|^2} - \frac{1}{\left|x-\frac{\lam^2 y }{|y|^2} \right|^2}\right] [u^{p_2}(y) e^{q_2 v(y)} - u_{\lam}^{p_2} e^{q_2 v_{\lam}} ] dy + \frac{2(\alpha-\mu)}{|x|^2}.
\end{align*}
\hfill $ \Box$

Now we show that the moving spheres method can be started at some $ \lam$.

\begin{pro}\label{b53}
Assume that $ \alpha > \mu$. Then for $ \lam$ large enough, we have
\[
\widetilde{S}^u_{\lam} (x) \geq 0, \quad \widetilde{S}^v_{\lam}(x) \geq 0
\]
and
\[
- \Delta \widetilde{S}^v_{\lam}(x) \geq 0
\]
for $ x \in B_{\lam}(0) \setminus \{0\}$.
\end{pro}

\noindent \emph{Proof.} Step 1: There exists an $ R_0 > 0$, such that
\[
\widetilde{S}^u_{\lam} (x) \geq 0, \quad \widetilde{S}^v_{\lam}(x) \geq 0  \quad  \text{and} \quad - \Delta \widetilde{S}^v_{\lam}(x) \geq 0
\]
for $ R_0 \leq |x| \leq \frac{\lam}{2} $ and $ \lam$ large enough.

By Lemma \ref{b48}, we find
\begin{align*}
\widetilde{S}^u_{\lam} (x) & =  u(x) -  \frac{\lam^2}{|x|^2} u \left(\frac{\lam^2 x}{|x|^2} \right) \\
& = \frac{\beta + o(1)}{|x|^2} - (\beta + o(1)) \frac{\lam^2}{|x|^2}\cdot \left( \frac{\lam^2}{|x|} \right)^{-2} \\
& = (\beta + o(1)) (|x|^{-2} - \lam^{-2}) \\
& \geq (\beta + o(1))\frac{3}{\lam^2} \\
& >0.
\end{align*}
By the asymptotic behaviour of $ v(x)$, we get
\begin{align*}
\widetilde{S}^v_{\lam} & = v(x) - v\left(\frac{\lam^2 x}{|x|^2}\right) + \mu \ln |x| - \mu \ln \lam \\
& =( -\alpha + o(1))\ln |x|  + (\alpha + o(1))(2\ln \lam - \ln |x|)  + \mu \ln |x| - \mu \ln \lam  \\
& = (2\alpha + o(1) - \mu )(\ln \lam - \ln |x|) \\
& > 0
\end{align*}
for $ \lam$ and $ R_0$ large enough.

Finally, by Lemma \ref{b51}, we obtain
\begin{align*}
- \Delta \widetilde{S}^v_{\lam} & = - \Delta v(x) + \frac{\lam^4}{|x|^4} \Delta v \left(\frac{\lam^2 x}{|x|^2} \right)  - 4\lam^2 \frac{x \cdot \nabla v \left(\frac{\lam^2 x}{|x|^2} \right)}{|x|^4} - \frac{2\mu}{|x|^2} \\
& = \frac{1}{|x|^2} \left[ - \Delta v(x) |x|^2 +  \frac{\lam^4}{|x|^2} \Delta v \left(\frac{\lam^2 x}{|x|^2} \right)  - 4\lam^2 \frac{x}{|x|^2} \cdot \nabla v \left(\frac{\lam^2 x}{|x|^2}\right) - 2 \mu \right] \\
& = \frac{1}{|x|^2} [ 4\alpha -2 \mu + o(1) ]\\
& > 0
\end{align*}
for $ \lam$ and $ R_0$ large enough.

Step 2:  $ \widetilde{S}_\lam^u(x) \geq 0$, $ \widetilde{S}_\lam^v(x) \geq 0$ and $ -\Delta \widetilde{S}_\lam^v(x) \geq 0$ for $ \frac{\lam}{2} \leq |x| \leq \lam$.

We can get $ -\Delta \widetilde{S}_\lam^v(x) \geq 0$ in the same way as Step 1. Since $ \widetilde{S}_\lam^v$ satisfies
\begin{equation}
\left\{
\begin{aligned}
&-\Delta \widetilde{S}_\lam^v \geq 0 \quad \text{in} \ B_\lam \setminus B_{\frac{\lam}{2}} \\
&\widetilde{S}_\lam^v = 0 \qquad \text{on} \ \ptl B_\lam \\
&\widetilde{S}_\lam^v \geq 0 \qquad \text{on} \ \ptl B_{\frac{\lam}{2}},
\nonumber
\end{aligned}
\right.
\end{equation}
then we infer from the maximum principle that $ \widetilde{S}_\lam^v \geq 0$ for $  \frac{\lam}{2} \leq |x| \leq \lam$.

By a direct calculation, we observe that
\begin{align*}
-\Delta \widetilde{S}_\lam^u(x) &  =  u^{p_1} e^{q_1 v} - u_{\lam}^{p_1} e^{q_1 v_{\lam}} \\
& = u^{p_1}( e^{q_1 v} - e^{q_1 v_{\lam}}) + (u^{p_1} - u_{\lam}^{p_1}) e^{q_1 v_{\lam}} \\
& = q_1 u^{p_1} e^{q_1 \xi_v}(v- v_{\lam}) + p_1 \xi_u^{p_1 - 1} (u-u_{\lam}) e^{q_1 v_{\lam}}
\end{align*}
for $ x \in \mathbb{R}^4 \setminus \{0\}$, where $ \xi_v$ is between $ v(x)$ and $ v_{\lam}(x)$ and $ \xi_u$ is between $ u(x)$ and $ u_{\lam}(x)$. We see that $ q_1 u^{p_1} e^{q_1 \xi_v}(v- v_{\lam}) \geq 0$ and $ p_1 \xi_u^{p_1 - 1}e^{q_1 v_{\lam}} \geq 0$ in $ B_\lam \setminus B_{\frac{\lam}{2}}$. Then
\begin{equation}
\left\{
\begin{aligned}
& - \Delta \widetilde{S}_\lam^u(x) - c(x) \widetilde{S}_\lam^u(x) \geq 0 \quad \text{in}  \ B_\lam \setminus B_{\frac{\lam}{2}} \\
&  \widetilde{S}_\lam^u(x) = 0 \qquad \text{on} \ \ptl B_\lam \\
&\widetilde{S}_\lam^u(x) \geq 0 \qquad \text{on} \ \ptl B_{\frac{\lam}{2}},
\nonumber
\end{aligned}
\right.
\end{equation}
where $ c(x)=p_1 \xi_u^{p_1 - 1}e^{q_1 v_{\lam}} > 0$ is bounded. By Lemma \ref{b48} and $ v(x) \sim -\alpha \ln |x|$ at $ \infty $, we find $ c(x) \leq \frac{C}{|x|^4}$, $ x \in  B_\lam \setminus B_{\frac{\lam}{2}} $. Similar argument as in Step 2 of Proposition \ref{b40}, one obtains $ \widetilde{S}_\lam^u(x) \geq 0$ for $  \frac{\lam}{2} \leq |x| \leq \lam$.

Step 3:  $ \widetilde{S}_\lam^u(x) \geq 0$, $ \widetilde{S}_\lam^v(x) \geq 0$ and $ -\Delta \widetilde{S}_\lam^v(x) \geq 0$ in $ B_{R_0} \setminus \{0\}$.

Since
\[
v\left(\frac{\lam^2 x}{|x|^2}\right)=[-2\alpha + o(1)]\ln \lam + [\alpha + o(1)] \ln |x|
\]
for $ \lam$ large enough, then
\[
\widetilde{S}_\lam^v(x) =v(x) + [2\alpha -\mu + o(1)]\ln \lam + [\mu -\alpha + o(1)]\ln |x|.
\]
Hence, we choose $ \lam$ large enough, such that $ \widetilde{S}_\lam^v(x) \geq 0$ for $ 0< |x| < R_0$.

Similarly, it follows from Lemma \ref{b48} that
\[
\widetilde{S}_\lam^u(x) = u(x) - \frac{\lam^2}{|x|^2}u\left(\frac{\lam^2 x}{|x|^2}\right) = u(x) - \frac{\beta + o(1)}{\lam^2}> 0
\]
in $ B_{R_0}\setminus \{0\}$ for $ \lam$ large enough.

By Lemma \ref{b52} (ii), we have $ - \Delta \widetilde{S}_\lam^v(x) \geq 0$ in $ B_{R_0}\setminus \{0\}$. \hfill $ \Box$

Now we define
\[
\widetilde{S}^u_{\lam,b}(x) = u_b(x) -  \frac{\lam^2}{|x|^2} u_b \left(\frac{\lam^2 x}{|x|^2} \right) \quad \text{and} \quad
\widetilde{S}^v_{\lam,b}(x) = v_b(x) - v_b\left(\frac{\lam^2 x}{|x|^2} \right) + \mu \ln|x| - \mu \ln \lam.
\]
For fixed $ b\in \mathbb{R}^4$, we define
\[
\lam_b = \inf\{\lam > 0 \ | \ \widetilde{S}^u_{\mu,b}(x)\geq 0,  \widetilde{S}^v_{\mu,b}(x)\geq 0, -\Delta  \widetilde{S}^v_{\mu,b}(x)\geq 0 \ \text{in} \ B_\mu \setminus \{0\}, \forall \lam\leq \mu < \infty\}.
\]

\begin{pro}\label{b54}
There exists a vector $ b \in \mathbb{R}^4$, such that $ \lam_b > 0$.
\end{pro}

\noindent \emph{Proof.} The proof of Proposition \ref{b54} is the same as the proof of Proposition \ref{b55}. We omit the details.  \hfill $ \Box$

\begin{pro}\label{b58}
If $ (u,v)$ is a nontrivial solution of problem \eqref{b1}, then $ \alpha > \mu$ can't occur.
\end{pro}

\noindent \emph{Proof.} Suppose on the contrary that $ \alpha > \mu$, then it follows from Proposition \ref{b54} that there exists a vector $ b \in \mathbb{R}^4$, such that $ \lam_b > 0$. Without loss of generality, one may assume $b=0$. Then it from the definition of $ \lam_0$ that
\be\label{b56}
\widetilde{S}_{\lam_0}^u(x) \geq 0 \quad \text{and} \quad \widetilde{S}_{\lam_0}^v(x)\geq 0 \quad \text{in} \quad B_{\lam_0} \setminus \{0\}.
\ee
Moreover, we deduce from Lemma \ref{b52} (ii) that
\[
- \Delta \widetilde{S}_{\lam_0}^v(x) > 0 \quad \text{in} \quad B_{\lam_0} \setminus \{0\}.
\]
Then the maximum principles implies
\[
\widetilde{S}_{\lam_0}^v(x) > 0 \quad \text{in} \quad B_{\lam_0} \setminus \{0\}.
\]
By \eqref{b56}, we know that
\[
- \Delta \widetilde{S}_{\lam_0}^u(x) = u^{p_1}(x) e^{q_1 v(x)} - u_{\lam_0}^{p_1}(x) e^{q_1 v_{\lam_0}(x)} \geq 0  \quad \text{in} \quad B_{\lam_0} \setminus \{0\}
\]
and
\[
\widetilde{S}_{\lam_0}^u(x) > 0 \quad \text{in} \quad B_{\lam_0} \setminus \{0\}.
\]

By the definition of $ \lam_0$, one of the following three cases may occur.

(i) $ \exists \lam_k < \lam_0$, $ \lam_k \rightarrow \lam_0$ with $ \inf\limits_{B_{\lam_k}\setminus \{0\}}  \widetilde{S}_{\lam_k}^u(x) < 0$.

(ii) $ \exists \lam_k < \lam_0$, $ \lam_k \rightarrow \lam_0$ with $ \inf\limits_{B_{\lam_k}\setminus \{0\}} \widetilde{S}_{\lam_k}^v(x) < 0$.

(iii) $ \exists \lam_k < \lam_0$, $ \lam_k \rightarrow \lam_0$, such that $ \widetilde{S}_{\lam_k}^u(x) \geq 0$, $ \widetilde{S}_{\lam_k}^v(x) \geq 0$ for $ x \in B_{\lam_k}\setminus \{0\}$, but $ \inf\limits_{B_{\lam_k}\setminus \{0\}} - \Delta \widetilde{S}_{\lam_k}^v(x) < 0$. \\
We will show that each case will lead to a contradiction.

If case (i) occur, then we deduce from $ \widetilde{S}_{\lam_0}^u(x) > 0 $ in $ B_{\lam_0} \setminus \{0\} $ and the Hopf Lemma that
\be\label{b57}
\frac{\ptl \widetilde{S}_{\lam_0}^u(x)}{\ptl \nu}  < 0 \quad \text{on} \quad \ptl B_{\lam_0},
\ee
where $ \nu$ is the unit outer normal direction. We define
\[
\gamma = \min_{\ptl B_{\frac{\lam_0}{2}}}\widetilde{S}_{\lam_0}^u(x) > 0
\]
and
\[
h(x) = \gamma - \frac{r^2}{|x|^2}\gamma \quad \text{in} \quad B_{\frac{\lam_0}{2}} \setminus B_r
\]
with $ r$ small. Then $ k(x) := \widetilde{S}_{\lam_0}^u(x) - h(x)$ satisfies
\begin{equation}
\left\{
\begin{aligned}
& - \Delta k(x) = - \Delta \widetilde{S}_{\lam_0}^u(x) \geq 0 \quad \text{in} \quad B_{\frac{\lam_0}{2}} \setminus B_r \\
& k(x) = \widetilde{S}_{\lam_0}^u(x) > 0 \qquad \text{on} \quad \ptl B_r \\
& k(x) > 0 \qquad \qquad \text{on} \quad \ptl B_{\frac{\lam_0}{2}}.
\nonumber
\end{aligned}
\right.
\end{equation}
Hence, we get from the maximum principle that $ \widetilde{S}_{\lam_0}^u(x) \geq \gamma - \frac{r^2}{|x|^2}\gamma$ in $ B_{\frac{\lam_0}{2}} \setminus B_r $. Letting $ r \rightarrow 0$, we get $ \widetilde{S}_{\lam_0}^u(x) \geq \gamma$ in $ B_{\frac{\lam_0}{2}} \setminus \{0\}$.
Hence, $ \inf\limits_{B_{\lam_k}\setminus \{0\}} \widetilde{S}_{\lam_k}^u(x) < 0$ is attained at some  $ x_k \in B_{\lam_k}\setminus B_{\frac{\lam_0}{2}} $ for $ k$ large enough. Moreover, one has
\[
\nabla \widetilde{S}_{\lam_k}^u(x_k) = 0.
\]
We can assume that, up to a subsequence, $ x_k \rightarrow x_0$, then $ \nabla \widetilde{S}_{\lam_0}^u(x_0)=0 $ and $ \widetilde{S}_{\lam_0}^u(x_0) =0$. Hence $ x_0 \in \ptl B_{\lam_0}$. However, this contradicts equation \eqref{b57}. Hence, case (i) can't occur.

Next, we show that case (ii) can not occur. Similarly, we define $ \hat{\gamma} = \min\limits_{\ptl B_{\frac{\lam_0}{2}}}\widetilde{S}_{\lam_0}^v(x) > 0$ and
\[
\hat{h}(x) = \hat{\gamma} - \frac{r^2}{|x|^2}\hat{\gamma} \quad \text{in} \quad B_{\frac{\lam_0}{2}} \setminus B_r
\]
with $ r$ small. Then $ \hat{k}(x) := \widetilde{S}_{\lam_0}^v(x) - \hat{h}(x)$ satisfies the condition of the maximum principle, we conclude that
\[
\widetilde{S}_{\lam_0}^v(x) > \hat{\gamma} > 0 \quad \text{in}  \quad B_{\frac{\lam_0}{2}} \setminus \{0\}.
\]

If case (ii) occur, then $ \inf\limits_{B_{\lam_k}\setminus \{0\}} \widetilde{S}^v_{\lam_k} (x) < 0$ is attained at some $ x_k \in B_{\lam_k} \setminus B_{\frac{\lam_0}{2}}$. Therefore,
\[
\nabla \widetilde{S}^v_{\lam_k}(x_k) = 0 .
\]
We can assume that, up to a subsequence, $ x_k \rightarrow x_0$ as $ k \rightarrow \infty$, then
\[
\nabla \widetilde{S}^v_{\lam_0}(x_0) = 0  \quad \text{and} \quad  \widetilde{S}_{\lam_0}^v(x_0) = 0,
\]
which implies that $ x_0 \in \ptl B_{\lam_0}$. But this contradicts the Hopf Lemma.

Finally, we show that case (iii) can't occur. We deduce from Lemma \ref{b52} (ii) and $ \widetilde{S}_{\lam_k}^u(x) \geq 0$, $ \widetilde{S}_{\lam_k}^v(x) \geq 0$ that
\[
- \Delta \widetilde{S}_{\lam_k}^v  > 0 \quad \text{in} \quad B_{\lam_k}\setminus \{0\},
\]
which contradicts
\[
\inf\limits_{B_{\lam_k}\setminus \{0\}} - \Delta \widetilde{S}_{\lam_k}^v(x) < 0.
\]
\hfill $ \Box$

\section{Proof of Theorem \ref{b59}}
By Proposition \ref{b46} and Proposition \ref{b58}, we have the following result.

\begin{pro}\label{b19}
Let $ (u,v)$ be a nontrivial solution for system \eqref{b1} and $ \alpha$ be defined as in \eqref{b3}. Then we have
\[
\alpha = \mu.
\]
\end{pro}

From Proposition \ref{b19}, we are now in the position to prove Theorem \ref{b59}. Most of the proof is similar as the proof in Section 4. Hence we will sketch most of the proof in the following.
\begin{pro}\label{b62}
Let $ \widetilde{S}^u_{\lam} (x)$ and $ \widetilde{S}^v_{\lam}(x)$ be defined as in Section 4. Then there exists a $ \lam > 0$ large enough, such that
\[
\widetilde{S}^u_{\lam} (x) \geq 0, \quad \widetilde{S}^v_{\lam}(x) \geq 0
 \quad \text{and} \quad
- \Delta \widetilde{S}^v_{\lam}(x) \geq 0 \quad  \text{in} \quad B_{\lam}(0) \setminus \{0\}.
\]
\end{pro}

\noindent \emph{Proof.} The proof of Proposition \ref{b62} is similar to the proof of Proposition \ref{b53}. \hfill $ \Box$

The definition of $ \lam_b$ is the same as Section 4.
\begin{pro}\label{b63}
Let $ \widetilde{S}^u_{\lam} (x)$ and $ \widetilde{S}^v_{\lam}(x)$ be defined as in Section 4. Then there exists a vector $ \bar{b} \in \mathbb{R}^4$, such that $ \lam_{\bar{b}} > 0$.
\end{pro}

\noindent \emph{Proof.} The proof of Proposition \ref{b63} is same as the proof of Proposition \ref{b54}. \hfill $ \Box$

\begin{pro}\label{b64}
Suppose that $ \lam_b > 0$ for some $ b \in \mathbb{R}^4$, then we have $
\widetilde{S}_{\lam_b,b}^u(x)\equiv 0$ and $ \widetilde{S}_{\lam_b,b}^v(x)\equiv 0$
in $  B_{\lam_b} \setminus \{0\}$.
\end{pro}

\noindent \emph{Proof.}
Without loss of generality, one may assume $b=0$. Then it from the definition of $ \lam_0$ that
\be
\widetilde{S}_{\lam_0}^u(x) \geq 0 \quad \text{and} \quad \widetilde{S}_{\lam_0}^v(x)\geq 0 \quad \text{in} \quad B_{\lam_0} \setminus \{0\}.
\ee
Suppose on the contrary that $ \widetilde{S}_{\lam_0}^u(x) \not\equiv 0$ or $ \widetilde{S}_{\lam_0}^v(x) \not\equiv  0$, then we deduce from Lemma \ref{b52} (ii) that
\[
- \Delta \widetilde{S}_{\lam_0}^v(x) > 0 \quad \text{in} \quad B_{\lam_0} \setminus \{0\}.
\]
The rest part of the proof is similar to that of Proposition \ref{b58}.  \hfill $ \Box$

\begin{pro}\label{b69}
For all $ b \in \mathbb{R}^4 $, we have $ \lam_b > 0$.
\end{pro}

\noindent \emph{Proof.} By Proposition \ref{b64}, there exists a vector $ \bar{b} \in \mathbb{R}^4 $ such that $ \lam_{\bar{b}} > 0$ and
\be\label{b65}
v_{\bar{b}}(x) - v_{\bar{b}}\left(\frac{\lam_{\bar{b}}^2 x}{|x|^2} \right) + \mu \ln|x| - \mu \ln \lam_{\bar{b}} \equiv 0 \quad \text{in} \quad B_{\lam_{\bar{b}}} \setminus \{0\}.
\ee
Letting $ |x|\rightarrow \infty$ in \eqref{b65}, then
\be\label{b68}
\lim\limits_{|x|\rightarrow \infty}( v_{\bar{b}}(x)+ \mu \ln|x|) = v_{\bar{b}}(0) + \mu \ln \lam_{\bar{b}}.
\ee
Suppose on the contrary that there exists a vector $ b \in \mathbb{R}^4$ such that $ \lam_b = 0$, then we obtain
\be\label{b66}
v_{b}(x) - v_{b}\left(\frac{\lam^2 x}{|x|^2} \right) + \mu \ln|x| - \mu \ln \lam \geq 0
\ee
for all $ \lam > 0$ and $ x\in B_\lam \setminus \{0\}$. Fix $ \lam$ and let $ |x| \rightarrow 0$ in \eqref{b66}, then
\[
\lim_{|x| \rightarrow 0} \left( v_{b}\left(\frac{\lam^2 x}{|x|^2} \right)  - \mu \ln|x| \right) \leq v_b(0) - \mu \ln \lam
\]
or
\be\label{b71}
\lim_{|x| \rightarrow \infty} ( v_{b}(x)  + \mu \ln|x|) \leq v_b(0) + \mu \ln \lam.
\ee
We infer from \eqref{b68} and \eqref{b71} that
\[
v_{\bar{b}}(0) + \mu \ln \lam_{\bar{b}} \leq v_b(0) + \mu \ln \lam,
\]
which is impossible as $\lambda\rightarrow 0$.
 \hfill $ \Box$

\begin{pro}\label{b70}
For all $ b \in \mathbb{R}^4 $, we have $ \lam_b > 0$ and $ \widetilde{S}^u_{\lam_b,b} \equiv 0$, $ \widetilde{S}^v_{\lam_b,b} \equiv 0$ in $ B_{\lam_b}\setminus \{0\}$.
\end{pro}
\noindent \emph{Proof.} This is a direct consequence of Proposition \ref{b64}
and Proposition \ref{b69}. \hfill $ \Box$

\noindent \emph{Proof of Theorem \ref{b59}.} Define $ f(x) = e^{v(x)}$, then it follows from Proposition \ref{b70} that
\[
f(x) = \frac{\lam_b^\mu}{|x-b|^\mu} f \left(\frac{\lam_b^2(x-b)}{|x-b|^2}+b\right)
\]
and
\[
u(x) = \frac{\lam_b^2}{|x-b|^2} u \left(\frac{\lam_b^2(x-b)}{|x-b|^2}+b\right).
\]
Let
\[
A= \lim_{|x|\rightarrow \infty} |x|^\mu f(x) = \lam_b^\mu f(b)
\]
and
\[
B = \lim_{|x|\rightarrow \infty} |x|^2 u(x) = \lam_b^2 u(b),
\]
where $ A$, $ B > 0$. We first assume that $ A=B=1$. Since
\[
u(x) =\frac{\lam_0^2}{|x|^2} u \left(\frac{\lam_0^2 x}{|x|^2}\right) = \frac{\lam_b^2}{|x-b|^2} u \left(\frac{\lam_b^2(x-b)}{|x-b|^2}+b\right),
\]
then one has
\be\label{b16}
u(x) =\frac{\lam_0^2}{|x|^2} \left[u(0)+ \nabla u(0) \frac{\lam_0^2 x}{|x|^2} + o\left(\frac{1}{|x|}\right)\right]
\ee
and
\be\label{b17}
u(x) =\frac{\lam_b^2}{|x-b|^2} \left[u(b)+ \nabla u(b) \frac{\lam_b^2(x-b)}{|x-b|^2} + o\left(\frac{1}{|x-b|}\right)\right]
\ee
as $ |x| \rightarrow \infty$. It follows from equation \eqref{b16} and \eqref{b17} that
\[
\frac{\ptl u(b)}{\ptl x_i} u(b)^{-2} = \frac{\ptl u(0)}{\ptl x_i} u(0)^{-2} - 2b_i
\]
and
\[
(u^{-1})_i(b) = 2b_i + (u^{-1})_i(0) = \frac{\ptl}{\ptl b_i}( |b|^2 + \nabla u^{-1}(0) \cdot b).
\]
Therefore,
\[
u(b) = \frac{1}{|b-d_0|^2 + d}.
\]
Similarly, we conclude that
\[
f(b)= \left( \frac{1}{|b-d_0|^2 + d} \right)^{\frac{\mu}{2}}
\]
and
\[
v(b) =\frac{\mu}{2} \ln \left( \frac{1}{|b-d_0|^2 + d} \right).
\]
Finally, if we don't assume $ A = B =1$, then
\[
u(x)= \frac{C_1(\eps)}{|x-x_0|^2 + \eps^2}
\]
and
\[
v(x)= \frac{3-p_1}{q_1}\ln \left( \frac{C_2(\eps)}{|x-x_0|^2 + \eps^2} \right).
\]
Then, by Proposition \ref{b19} and direct calculations, one obtains
\begin{equation}
\left\{
\begin{aligned}
& C_1^{p_1} C_2^{3-p_1} = 8C_1 \eps^2 \\
& C_1^{p_2} C_2^{4-p_2} = 48 \mu \eps^4.
\nonumber
\end{aligned}
\right.
\end{equation}

This completes the proof of Theorem \ref{b59}. \hfill $ \Box$

\section{Acknowledgments.}
 This work is partially supported by by National Natural Science Foundation of China 11871160 and 12141105.

\bibliographystyle{plain}
\bibliography{citezu}

\end{document}